# Permanents, Pfaffian orientations, and even directed circuits

By Neil Robertson,* P. D. Seymour, and Robin Thomas**

**Abstract**

Given a 0-1 square matrix $A$, when can some of the 1's be changed to $-1$'s in such a way that the permanent of $A$ equals the determinant of the modified matrix? When does a real square matrix have the property that every real matrix with the same sign pattern (that is, the corresponding entries either have the same sign or are both zero) is nonsingular? When is a hypergraph with $n$ vertices and $n$ hyperedges minimally nonbipartite? When does a bipartite graph have a "Pfaffian orientation"? Given a digraph, does it have no directed circuit of even length? Given a digraph, does it have a subdivision with no even directed circuit?

It is known that all of the above problems are equivalent. We prove a structural characterization of the feasible instances, which implies a polynomial-time algorithm to solve all of the above problems. The structural characterization says, roughly speaking, that a bipartite graph has a Pfaffian orientation if and only if it can be obtained by piecing together (in a specified way) planar bipartite graphs and one sporadic nonplanar bipartite graph.

## 1. Introduction

Computing the permanent of a matrix seems to be of a different computational complexity from computing the determinant. While the determinant can be calculated using Gaussian elimination, no efficient algorithm for computing the permanent is known, and, in fact, none is believed to exist. More precisely, Valiant [24] has shown that computing the permanent is #P-complete even when restricted to 0-1 matrices.

*Research partially performed under a consulting agreement with Bellcore, and partially supported by DIMACS Center, Rutgers University, New Brunswick, New Jersey 08903, USA, by NSF under Grant No. DMS-9401981, and by ONR under Grant No. N00014-92-J-1965.

**Research partially performed under a consulting agreement with Bellcore, and partially supported by DIMACS Center, Rutgers University, New Brunswick, New Jersey 08903, USA, by NSF under Grant No. DMS-9623031, and by ONR under Grant No. N00014-93-1-0325.

It is therefore reasonable to ask if perhaps computing the permanent can be somehow reduced to computing the determinant of a related matrix. In particular, the following question was asked by Pólya [17] in 1913. If $A$ is a 0-1 square matrix, under what conditions does there exist a matrix $B$ obtained from $A$ by changing some of the 1's to $-1$'s in such a way that the permanent of $A$ equals the determinant of $B$? For the purpose of this paper let us say that $B$ (when it exists) is a *Pólya matrix* for $A$. The complexity status of the decision problem of whether an input matrix has a Pólya matrix remained open until present. In this paper we solve the problem. Specifically, we first give a structural characterization of matrices that have a Pólya matrix. Roughly speaking, they can all be obtained by piecing together "planar" matrices and one sporadic nonplanar matrix. We then use this characterization to design a polynomial-time algorithm that given an input matrix $A$ outputs either a Pólya matrix for $A$, or a certain "obstruction" submatrix of $A$ whose presence implies that $A$ has no Pólya matrix. The algorithm easily extends to matrices with nonnegative entries, as pointed out by Vazirani and Yannakakis [25].

Our results are best stated and proved in terms of bipartite graphs. By a *graph* we mean a finite simple undirected graph, that is, one with no loops or parallel edges. A set $M$ of edges of $G$ is a *matching* if every vertex of $G$ is incident with at most one edge in $M$; it is a *perfect matching* if every vertex of $G$ is incident with exactly one edge in $M$. Let $G$ be a graph, and let $H$ be a subgraph of $G$. We say that $H$ is *central* if $G\setminus V(H)$ has a perfect matching. (If $G$ is a graph, and $X$ is a vertex, an edge, or a set of vertices or edges, then $G\setminus X$ denotes the graph obtained from $G$ by deleting $X$.) Let $D$ be an orientation of $G$, and let $C$ be a circuit of $G$ of even length. (*Paths* and *circuits* have no "repeated" vertices.) We say that $C$ is *oddly oriented* (in $D$) if $C$ contains an odd number of edges that are directed (in $D$) in the direction of each orientation of $C$. We say that $D$ is a *Pfaffian orientation* of $G$ if every central circuit of $G$ of even length is oddly oriented in $D$.

A graph $G$ is *bipartite* if its vertex-set can be partitioned into two sets $A$ and $B$ in such a way that every edge has one end in $A$ and the other end in $B$. We say that $(A, B)$ is a *bipartition* of $G$, and we refer to $A$ and $B$ as *color classes*. With every 0-1 square matrix $A$ we associate a bipartite graph $G$ as follows. There is a vertex of $G$ corresponding to every row and every column of $A$, and two vertices of $G$ are adjacent if and only if one represents a row, say $r$, and the other represents a column, say $c$, such that the entry of $A$ in row $r$ and column $c$ is nonzero. Vazirani and Yannakakis [25] proved the following.

1.1. *Let $A$ be a 0-1 square matrix, and let $G$ be the associated bipartite graph. Then $A$ has a Pólya matrix if and only if $G$ has a Pfaffian orientation.*

Little [9] proved the following elegant characterization of bipartite graphs that admit a Pfaffian orientation (and hence of matrices that admit a Pólya matrix). We say that a graph $G$ is a *subdivision* of a graph $H$ if $G$ is obtained from $H$ by replacing the edges of $H$ by internally disjoint paths, each containing at least one edge. We say that $G$ is an *even subdivision* of $H$ if $G$ is obtained from $H$ by replacing the edges of $H$ by internally disjoint paths, each containing an even number of vertices and at least one edge. We say that a graph $G$ *contains* a graph $H$ and that $H$ is *contained* in $G$ if some even subdivision of $H$ is isomorphic to a central subgraph of $G$.

1.2. *A bipartite graph admits a Pfaffian orientation if and only if it does not contain $K_{3,3}$.*

Unfortunately, 1.2 does not seem to imply a polynomial-time algorithm to test whether a bipartite graph has a Pfaffian orientation, the difficulty being that it is not clear how to efficiently test for a $K_{3,3}$ containment. Our main result gives a structural description of graphs that admit a Pfaffian orientation, and it enabled us to derive a polynomial-time recognition algorithm.

To state our main result we need some definitions. Let $G_0$ be a graph, let $C$ be a central circuit of $G_0$ of length four, and let $G_1, G_2$ be two subgraphs of $G_0$ such that $G_1 \cup G_2 = G_0$, $G_1 \cap G_2 = C$, $V(G_1) - V(G_2) \neq \emptyset$ and $V(G_2) - V(G_1) \neq \emptyset$. (The intersection and union of two subgraphs of a graph is defined in the natural way.) Let $G$ be obtained from $G_0$ by deleting some (possibly none) of the edges of $C$. In these circumstances we say that $G$ is a *4-sum* of $G_1$ and $G_2$. The *Heawood graph* is the bipartite graph associated with the incidence matrix of the Fano plane (see Figure 1 below).

A graph $G$ is *k-extendable*, where $k \geq 0$ is an integer, if every matching of size at most $k$ can be extended to a perfect matching. A 2-extendable bipartite graph is called a *brace*. It is easy to see (and will be outlined in §7) that the problem of finding Pfaffian orientations of bipartite graphs can be reduced to braces. The following is our main result.

1.3. *A brace has a Pfaffian orientation if and only if either it is isomorphic to the Heawood graph, or it can be obtained from planar braces by repeated application of the 4-sum operation.*

In Section 9 we use 1.3 to design a polynomial-time algorithm to decide if a bipartite graph has a Pfaffian orientation. By 1.1 this also gives a polynomial-time algorithm to decide if a 1-1 square matrix has a Pólya matrix, by [25] this solves the even circuit problem for directed graphs (see [19], [20], [21], [22], [23]), by [11], [18] it solves the problem of determining which hypergraphs with $n$ vertices and $n$ hyperedges are minimally nonbipartite, and by [8] it solves the problem of determining which real $n \times n$ matrices are sign-nonsingular.

932    NEIL ROBERTSON, P. D. SEYMOUR, AND ROBIN THOMAS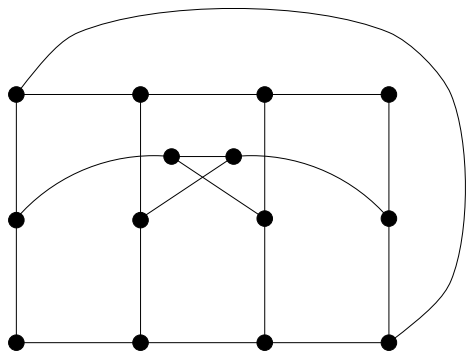

The Heawood graph

Figure 1

See also [1] for variations of sign-singularity. We briefly review some of these equivalent formulations in Section 7.

Let us outline the proof now. We first characterize containment-minimal nonplanar 1-extendable bipartite graphs. To state this characterization we need to define several classes of graphs. Let $k \geq 2$ be an integer, and let $C$ be a circuit with vertices $u_1, u_2, \ldots, u_{2k}$ listed in their order on $C$. Let $G$ be the graph obtained from $C$ by adding $k$ edges $e_1, e_2, \ldots, e_k$ where $e_i$ has ends $u_i$ and $u_{i+k}$ ($i = 1, 2, \ldots, k$). If $k \geq 4$ we say that $G$ is a *Möbius ladder*; the edges $e_1, e_2, \ldots, e_k$ are called *rungs*. Let $H$ be the graph obtained from $C$ by adding two vertices $v_1$ and $v_2$, an edge joining $v_1$ and $v_2$, and for $i = 1, 2, \ldots, k$ an edge with ends $v_1$ and $u_{2i-1}$, and an edge with ends $v_2$ and $u_{2i}$. We say that $H$ is a *biwheel*; the vertices $v_1$ and $v_2$ are called *hubs*, and the edges incident with exactly one of them are called *spokes*. A *stem* is a graph obtained from a Möbius ladder with an even number of rungs by replacing all rungs by disjoint paths on three vertices with the same ends. A *flower* is a graph obtained from a biwheel by replacing all spokes by internally disjoint paths on three vertices with the same ends. The flower on ten vertices is called *Bud*. Let $(\{a_1, a_2, a_3\}, \{b_1, b_2, b_3\})$ be the bipartition of $K_{3,3}$. We define Uno to be the graph obtained from $K_{3,3}$ by subdividing every edge incident with $a_1$ exactly once, and adding a two-edge path joining $a_2$ and $a_3$. We define Duo to be the graph obtained from $K_{3,3}$ by adding a two-edge path joining $a_2$ and $a_3$, and a two-edge path joining $b_2$ and $b_3$. See Figure 2.

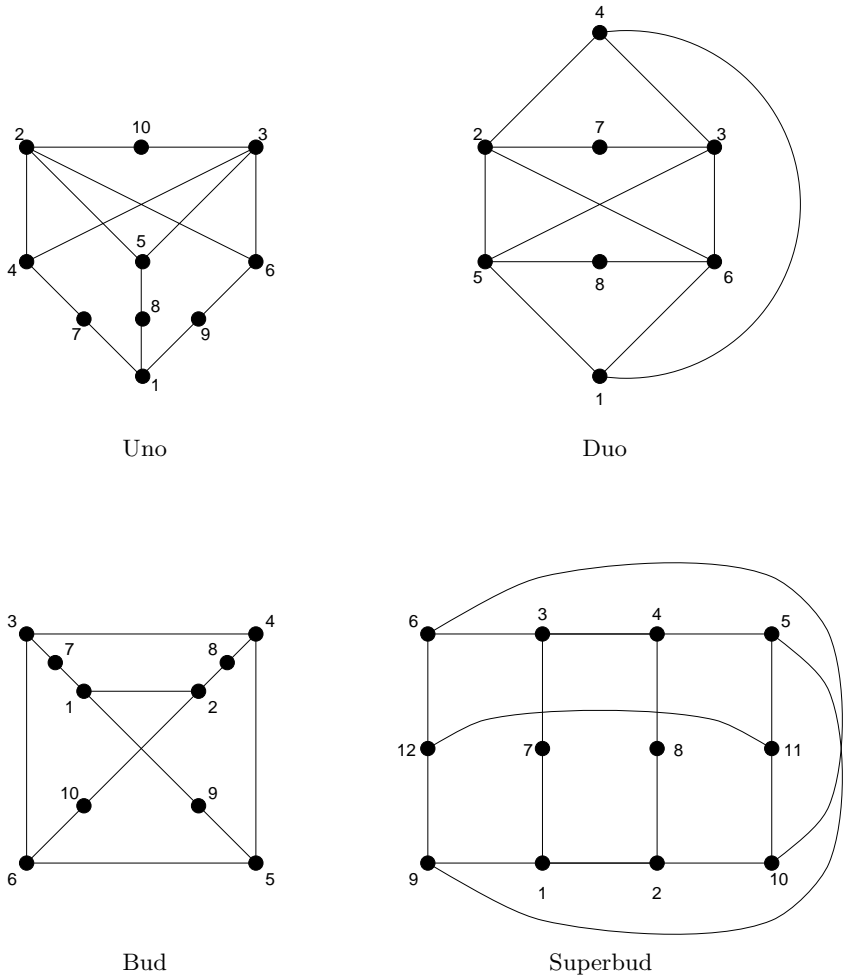

Figure 2

1.4. *Let G be a nonplanar 1-extendable bipartite graph. Then G contains one of the following graphs*:

(i) $K_{3,3}$,

(ii) *a stem*,

(iii) *a flower*,

(iv) *Uno, or*

(v) *Duo.*

To prove 1.4 we first give in Section 2 a description of graphs with no "magic" circuits. We then prove in Section 3 that if $G$ is as in 1.4 and is containment-minimal, then it has no magic circuit. Using the result of Section 2 it is then easy to prove 1.4.

As a next step we use 1.4 in Section 4 to prove the following.

1.5. *Let $G$ be a nonplanar 2-extendable bipartite graph. Then $G$ contains one of the following graphs:*
(i) *$K_{3,3}$,*
(ii) *the Heawood graph, or*
(iii) *Rotunda.*

The graph Rotunda is defined as follows. Let $C$ be a circuit of length four, and let $H$ be obtained from $C$ by adding four new vertices of degree one, each adjacent to a different vertex of $C$. Let the new vertices be $a, b, c, d$ listed in the order of their neighbors on $C$. Let $H_1, H_2, H_3$ be three isomorphic copies of $H$, and let $a_i, b_i, c_i, d_i$ ($i = 1, 2, 3$) be the vertices corresponding to $a, b, c, d$, respectively. Then Rotunda is obtained by identifying $a_1, a_2, a_3$ into $a_0$, identifying $b_1, b_2, b_3$ into $b_0$, and so on. We say that $\{a_0, b_0, c_0, d_0\}$ is its *center*. See Figure 3 below.

In Section 5 we prove that if a brace $G$ does not contain $K_{3,3}$, but contains Rotunda, and $X$ is the set of vertices of $G$ corresponding to the center of Rotunda, then $G\backslash X$ is disconnected. This is a lemma that will be used in the next section to show that in those circumstances $G$ is a 4-sum of two smaller graphs. In Section 6 we derive our main theorem from 1.5, as follows. The easy "if" part follows from 1.2. To prove the difficult "only if" part of 1.3 we may assume that $G$ is a nonplanar brace with a Pfaffian orientation. Then by 1.2 and 1.5 $G$ contains the Heawood graph or Rotunda. We prove that if $G$ contains the Heawood graph, then it is isomorphic to it, and that if it contains Rotunda, then it is a 4-sum of two smaller braces. For the latter assertion we use the result of Section 5. In Section 7 we deduce several consequences of the main theorem, in Section 8 we prove several lemmas needed for the algorithm, and in Section 9 we design a polynomial-time algorithm to test if a bipartite graph has a Pfaffian orientation.

## 2. Magic circuits

An edge $e$ of a graph $G$ is called *planarizing* (in $G$) if $G\backslash e$ is planar, and *nonplanarizing* otherwise. We say that a circuit $C$ in a graph $G$ is *magic* if at most one edge of $C$ has the property that it is planarizing and has both ends of degree at least three. Our objective in this section is to prove a theorem about 2-connected nonplanar graphs with no magic circuits.

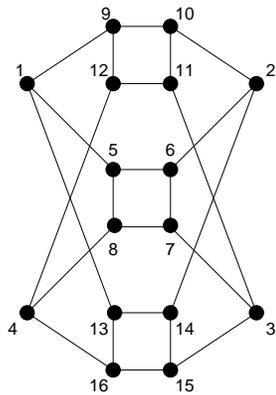

Rotunda

Figure 3

Let $(A, B)$ be the bipartition of $K_{3,3}$. Let $K_{3,3}^+$ denote the graph obtained from $K_{3,3}$ by joining two vertices of $A$ by an edge, let $K_{3,3}^{+^+}$ denote the graph obtained from $K_{3,3}^+$ by joining two vertices of $B$ by an edge, and let $K_{3,3}^{++}$ denote the graph obtained from $K_{3,3}$ by adding two edges joining two distinct pairs of vertices of $A$. The main result of this section is the following.

2.1. *Let $G$ be a 2-connected nonplanar graph with no magic circuit. Then $G$ is isomorphic to a subdivision of one of the following graphs*:
 (i) *a Möbius ladder*,
 (ii) *a biwheel*,
(iii) $K_{3,3}^+$,
(iv) $K_{3,3}^{+^+}$,
 (v) $K_{3,3}^{++}$, *or*
(vi) $K_5$.

Let $H$ be a subgraph of a graph $G$. By an $H$-*path* in $G$ we mean a path $P$ in $G$ with at least one edge such that $P$ has both ends in $H$ and is otherwise disjoint from $H$. A first step in the proof of 2.1 is the following. Let us remark that, by definition, every Möbius ladder has at least eight vertices.

2.2. *Let $G$ be a 2-connected graph with no magic circuit such that some subgraph of $G$ is isomorphic to a subdivision of a Möbius ladder. Then $G$ is isomorphic to a subdivision of a Möbius ladder.*

*Proof.* Choose a subgraph $H$ of $G$ such that

(i) $H$ is isomorphic to a subdivision of a Möobius ladder $L$, and

(ii) subject to (i), $|V(L)|$ is maximum.

Let $V(L) = \{u_1, u_2, \ldots, u_{2k}\}$ as in the definition of a Möbius ladder, let $v_1, v_2, \ldots, v_{2k}$ be the corresponding vertices of $H$, and for $i, j \in \{1, 2, \ldots, 2k\}$ such that $u_i$ and $u_j$ are adjacent in $L$ let $Q_{i,j}$ denote the corresponding path in $H$. We claim the following.

(1) *For every $H$-path $P$ in $G$ there exists an integer $i$ such that one of the ends of $P$ is equal to $v_i$ and the other belongs to $Q_{i+k-1,i+k} \cup Q_{i+k,i+k+1} \setminus v_{i+k}$ (index arithmetic is modulo $2k$).*

To prove (1) let $P$ be an $H$-path in $G$ with ends $u$ and $v$. We may assume that for some integer $i$ with $1 \leq i \leq k+1$ either $u \in V(Q_{1,k+1}) - \{v_1, v_{k+1}\}$, $v \in V(Q_{i,k+i}) \cup V(Q_{i,i+1}) - \{v_{i+1}, v_{k+i}\}$ and $i \leq k-1$, or $u \in V(Q_{1,2}) - \{v_2\}$ and $v \in V(Q_{i,i+1}) - \{v_{i+1}\}$. If $i < k-1$, then every edge of the circuit of $Q_{1,k+1} \cup Q_{1,2} \cup Q_{2,3} \cup \ldots \cup Q_{i,i+1} \cup Q_{i,k+i} \cup P$ is nonplanarizing, and hence the circuit is magic, a contradiction. If $i = k-1$, then every edge of the circuit $Q_{2,k+2} \cup Q_{2,3} \cup Q_{3,k+3} \cup Q_{k+2,k+3}$ is nonplanarizing. Thus we may assume that $i \in \{k, k+1\}$, $u \in V(Q_{1,2}) - \{v_2\}$ and $v \in V(Q_{i,i+1}) - \{v_{i+1}\}$. If $u = v_1$, then $v \neq v_{k+1}$, because otherwise $Q_{1,k+1} \cup P$ is magic, and so (1) holds. Thus we may assume that $u \in V(Q_{1,2}) - \{v_1, v_2\}$. Now $v \notin V(Q_{k+1,k+2}) - \{v_{k+1}, v_{k+2}\}$, because otherwise $H \cup P$ contradicts (ii). If $v = v_{k+1}$ then (1) holds, and so we may assume that $v \in V(Q_{k,k+1}) - \{v_{k+1}\}$. Now every edge of the circuit $Q_{2,3} \cup Q_{3,k+3} \cup Q_{k+2,k+3} \cup Q_{2,k+2}$ is nonplanarizing, and hence the circuit is magic, a contradiction. This proves (1).

We claim that $H = G$. To prove this suppose to the contrary that $H \neq G$; then (since $G$ is 2-connected) $G$ has an $H$-path $P$. By (1) we may assume that one end of $P$ is $v_1$ and the other is $v \in V(Q_{k+1,k+2}) - \{v_{k+1}\}$. Let $Q$ denote the subpath of $Q_{k+1,k+2}$ between $v_{k+1}$ and $v$. We may assume that $P$ is chosen so that $|V(Q)|$ is minimum; in other words there is no $H$-path $P'$ with ends $v_1$ and $v' \in V(Q) - \{v_{k+1}, v\}$. Since every edge of $Q_{1,k+1} \cup P$ is nonplanarizing and the circuit $Q_{1,k+1} \cup Q \cup P$ is not magic, we deduce that $V(Q) - \{v_{k+1}, v\}$ contains a vertex $w$ of degree at least three in $G$. Thus $w$ is an end of an $H$-path $P'$. By (1) and the choice of $P$ it follows that the other end of $P'$ is $v_2$ and that $P$ and $P'$ are vertex-disjoint. It is now easy to check that every edge of the circuit $C = Q_{3,4} \cup Q_{4,4+k} \cup Q_{3+k,4+k} \cup Q_{3,3+k}$ is nonplanarizing, and hence $C$ is magic, a contradiction. □

2.3. *Let $G$ be a 2-connected graph with no magic circuit, such that some subgraph of $G$ is isomorphic to a subdivision of a biwheel on at least eight vertices. Then $G$ is isomorphic to a subdivision of a biwheel.*

*Proof.* Choose a subgraph $H$ of $G$ such that

(i) $H$ is isomorphic to a subdivision of a biwheel $W$, and
(ii) subject to (i), $|V(W)|$ is maximum.

Let $u_0, u_{2k+1}$ be the hubs of $W$, let $u_1, u_2, \ldots, u_{2k}$ (in order) be the vertices of the circuit $W \setminus \{u_0, u_{2k+1}\}$ and assume that $u_0$ is adjacent to $u_1, u_3, \ldots, u_{2k-1}$. Thus $u_{2k+1}$ is adjacent to $u_2, u_4, \ldots, u_{2k}$. Let $v_0, v_1, v_2, \ldots, v_{2k}, v_{2k+1}$ be the corresponding vertices of $H$, and let $Q_{i,j}$ be the corresponding paths. Let $C$ denote the circuit $Q_{1,2} \cup Q_{2,3} \cup \cdots \cup Q_{2k-1,2k} \cup Q_{2k,1}$. We need the following claim.

(1) *Every $H$-path either has one end $v_0$ and the other in*
$$V(C) - \{v_1, v_3, \ldots, v_{2k-1}\},$$
*or it has one end $v_{2k+1}$ and the other in $V(C) - \{v_2, v_4, \ldots, v_{2k}\}$.*

To prove (1) let $P$ be an $H$-path with ends $w_1$ and $w_2$. Suppose first that $w_1, w_2 \notin \{v_0, v_{2k+1}\}$. Let $i \in \{1, 2\}$. If $w_i \in V(Q_{0,j})$ for some $j \in \{1, 3, \ldots, 2k-1\}$ or $w_i \in V(Q_{2k+1,j})$ for some $j \in \{2, 4, \ldots, 2k\}$ let $W_i$ be the subpath of $Q_{0,j}$ or $Q_{2k+1,j}$ with one end $w_i$ and the other end on $C$; otherwise let $W_i$ be the null graph. It follows that $W_1 \cup W_2 \cup C \cup P$ includes a circuit every edge of which is nonplanarizing, and hence the circuit is magic, a contradiction. Thus we may assume that $w_1 = v_0$. If $w_2 \in V(Q_{0,j})$ for some $j \in \{1, 3, \ldots, 2k-1\}$ then $Q_{0,j} \cup P$ includes a magic circuit, and if $w_2 \in V(Q_{2k+1,j}) - \{v_j\}$ for some $j \in \{2, 4, \ldots, 2k\}$, then $Q_{0,2k+1} \cup Q_{2k+1,j} \cup P$ includes a magic circuit, a contradiction. Thus $w_2 \in V(C) - \{v_1, v_3, \ldots, v_{2k-1}\}$, as desired. This proves (1).

We claim that $H = G$. To prove this suppose to the contrary that $H \neq G$; then (since $G$ is 2-connected) $G$ has an $H$-path $P$. By (1) we may assume that one end of $P$ is $v_0$, and that the other is $v \in V(C) - \{v_1, v_3, \ldots, v_{2k-1}\}$. We may assume that $v \in V(Q_{1,2}) - \{v_1\}$. Let $Q$ be the subpath of $Q_{1,2}$ between $v_1$ and $v$. We may assume that $P$ is chosen so that $|V(Q)|$ is minimum. In other words, there is no $H$-path $P'$ with one end $v_0$ and the other in $V(Q) - \{v\}$, and if $v = v_2$ then there is no $H$-path $P''$ with one end $v_{2k+1}$ and the other in $V(Q) - \{v_1, v_2\}$. Since the circuit $P \cup Q \cup Q_{0,1}$ is not magic, the path $Q$ has an internal vertex $w$ of degree at least three. Thus $w$ is an end of an $H$-path $P'$. By the minimality of $|V(Q)|$ the paths $P$ and $P'$ are disjoint. By (1) the other end of $P'$ is $v_0$ or $v_{2k+1}$; by the minimality of $Q$ the other end is $v_{2k+1}$ and $v \neq v_2$. Now $H \cup P \cup P'$ is isomorphic to a subdivision of the biwheel on $|V(W)| + 2$ vertices, contrary to (ii). This proves our claim that $G = H$, and hence completes the proof of the lemma. □

2.4. *Let $G$ be a 2-connected nonplanar graph with no magic circuit, such that no subgraph of $G$ is isomorphic to a subdivision of $K_{3,3}$. Then $G$ is isomorphic to a subdivision of $K_5$.*

*Proof.* By Kuratowski's theorem $G$ has a subgraph $K$ isomorphic to a subdivision of $K_5$. Let $v_1, v_2, \ldots, v_5$ be the corresponding vertices of $K$, and let $Q_{12}, Q_{13}, \ldots, Q_{45}$ be the corresponding paths. We claim that $G = K$. To prove this suppose to the contrary that $G \neq K$; then $G$ has a $K$-path $P$. If both ends of $P$ belong to $V(Q_{ij})$ for some distinct integers $i, j \in \{1, 2, \ldots, 5\}$, then $P$ and an appropriate subpath of $Q_{ij}$ form a circuit $C$ such that every edge of $C$ is nonplanarizing, and hence $C$ is magic, a contradiction. Thus the ends of $P$ belong to $V(Q_{ij})$ for no pair of distinct integers $i, j \in \{1, 2, \ldots, 5\}$, in which case it is easy to see that $K \cup P$ has a subgraph isomorphic to a subdivision of $K_{3,3}$, contrary to hypothesis. □

2.5. *Let $G$ be a 2-connected graph with no magic circuit, and let $K$ be a subgraph of $G$ isomorphic to a subdivision of $K_{3,3}$. Assume that no subgraph of $G$ is isomorphic to a Möbius ladder or a biwheel on at least eight vertices. Then $G$ is isomorphic to a subdivision of $K_{3,3}$, $K_{3,3}^+$, $K_{3,3}^{\overset{+}{+}}$, or $K_{3,3}^{++}$.*

*Proof.* Let $v_1, v_2, v_3, v_4, v_5, v_6$ be the vertices of $K$ of degree three in $K$, and let $Q_{ij}$ ($i = 1, 2, 3$, $j = 4, 5, 6$) be the corresponding paths, say. We claim the following.

(1) *For all $i = 1, 2, 3$, all $j = 4, 5, 6$ and every $K$-path $P$, not both ends of $P$ belong to $V(Q_{i4}) \cup V(Q_{i5}) \cup V(Q_{i6}) - \{v_4, v_5, v_6\}$ or $V(Q_{1j}) \cup V(Q_{2j}) \cup V(Q_{3j}) - \{v_1, v_2, v_3\}$.*

To prove (1) suppose for instance that a $K$-path $P$ has both ends in $V(Q_{i4}) \cup V(Q_{i5}) - \{v_4, v_5\}$. Then $Q_{i4} \cup Q_{i5} \cup P$ has a circuit $C$ such that every edge of $C$ is nonplanarizing in $G$. Thus $C$ is magic, a contradiction. This proves (1).

From (1) we deduce the following claim.

(2) *Every $K$-path $P$ has at least one end in $\{v_1, v_2, \ldots, v_6\}$.*

Indeed, otherwise $K \cup P$ is isomorphic to a subdivision of the Möbius ladder on eight vertices by (1), contrary to hypothesis. This proves (2).

To complete the proof assume first that every $K$-path has both ends in $\{v_1, v_2, \ldots, v_6\}$. Then, by (1), every $K$-path has both ends either in $\{v_1, v_2, v_3\}$ or in $\{v_4, v_5, v_6\}$. Furthermore, every two distinct $K$-paths are vertex-disjoint, except possibly for their ends, because otherwise $G$ would have a magic circuit. Similarly, no two distinct $K$-paths have the same set of ends. If $G$ has three $K$-paths with ends $v_1$ and $v_2$, $v_2$ and $v_3$, and $v_1$ and and $v_3$, then their union is a magic circuit, a contradiction. If $G$ has three $K$-paths $P_1, P_2, P_3$ with ends $v_1$ and $v_2$, $v_2$ and $v_3$, $v_4$ and $v_5$, respectively, then $P_1 \cup Q_{26} \cup Q_{16}$ is a magic circuit, a contradiction. (To see that $P_1 \cup Q_{26} \cup Q_{16}$ is magic notice that every edge of $P_1 \cup Q_{26}$ is nonplanarizing and that every internal vertex of $Q_{16}$ has degree two in $G$, because every $K$-path has both ends in $\{v_1, v_2, \ldots, v_6\}$.) It follows from the symmetry that $G$ is isomorphic to a subdivision of one of $K_{3,3}$, $K_{3,3}^+$, $K_{3,3}^{+\,+}$, $K_{3,3}^{++}$, as desired.

We may therefore assume that there is a $K$-path $P$ with one end $v_1$ and the other end $v \in V(Q_{24}) - \{v_2, v_4\}$. Let $C$ be the circuit of $Q_{14} \cup P \cup Q_{24}$. We may assume that $v$ and $P$ are chosen so that $|V(C) \cap V(Q_{24})|$ is minimum. Since every edge of $Q_{14} \cup P$ is nonplanarizing in $G$ and $C$ is not magic, we deduce that the subpath of $Q_{24}$ between $v_4$ and $v$ has an internal vertex $w$ of degree at least three. Thus $w$ is an end of a $K$-path $P'$. By (1) and (2) the other end of $P'$ is $v_k$ for some $k \in \{1, 3, 5, 6\}$. From the choice of $P$ we deduce that $P'$ is disjoint from $P \backslash v_1$. Now $k \neq 1$ by the choice of $v$, and $k \neq 3$, because otherwise $Q_{24} \cup Q_{34} \cup P'$ contains a circuit $C'$ such that every edge of $C'$ is nonplanarizing in $G$, in which case $C'$ is magic. Thus $k = 5$ or $k = 6$, and it follows that $K \cup P \cup P'$ is isomorphic to a subdivision of the biwheel on eight vertices, a contradiction. □

Theorem 2.1 follows immediately from 2.2, 2.3, 2.4 and 2.5.

## 3. One-extendable bipartite graphs

In this section we prove 1.4. If $G$ is a graph and $X \subseteq V(G)$, we denote by $N_G(X)$ the set of vertices of $V(G) - X$ adjacent to a vertex in $X$. We need the following well-known characterization of $k$-extendable bipartite graphs (see [13]).

3.1. *Let $G$ be a connected bipartite graph with bipartition $(A, B)$, and let $k \geq 0$ be an integer. Then the following two conditions are equivalent.*
(i) *$G$ is $k$-extendable, and*
(ii) *$|A| = |B|$, and for every nonempty subset $X$ of $A$, either $|N_G(X)| \geq |X| + k$, or $N_G(X) = B$ (or both).*

3.2. *Every 1-extendable connected bipartite graph with more than two vertices is 2-connected. Every connected brace on at least five vertices is 3-connected.*

*Proof.* This follows immediately from 3.1. □

The following lemma is crucial for the proof of 1.4. If $M$ is a matching in a graph $G$ and $P$ is a path or a circuit, we say that $P$ is *M-alternating* if every vertex of $P$ of degree two in $P$ is incident with an edge of $M \cap E(P)$.

3.3. *Let $G$ be a 1-extendable bipartite graph, let $M$ be a perfect matching of $G$, and let $C$ be a circuit in $G$ such that no two edges of $C$ form an edge cut of $G$. Then there exists an edge $e \in E(C) - M$ such that $G \backslash e$ is 1-extendable.*

*Proof.* We may assume that $G$ is connected. For $e \in E(G)$ let $a(e)$ denote the number of $M$-alternating circuits of $G$ that contain $e$. Choose $e \in E(C) - M$ with $a(e)$ minimum. We claim that $e$ satisfies the conclusion of the lemma. Suppose for a contradiction that it does not, and let $(A, B)$ be a bipartition of $G$. The graph $G \backslash e$ is connected by 3.2, and hence, by 3.1 applied to $G \backslash e$, there exists a nonempty proper subset $X$ of $A$ such that $|Y| = |X|$, where $Y = N_{G \backslash e}(X)$. Let $F$ be the set of all edges of $G$ with one end in $Y$ and the other in $A - X$. Then $F \cap M = \emptyset$, every $M$-alternating circuit containing a member of $F$ contains precisely one such member, and also contains $e$. Thus $a(e) \geq a(f)$ for every $f \in F$, and the inequality is strict if $|F| > 1$, because $a(f) > 0$ since $G$ is 1-extendable. Since $E(C) \cap F \neq \emptyset$ we deduce from the choice of $e$ that $|F| = 1$, say $F = \{f\}$. Then $\{e, f\}$ is an edge cut of $G$, contrary to the hypothesis of the lemma. Thus $e$ satisfies the conclusion of the lemma. □

Let $\mathcal{G}$ denote the class of all containment-minimal 1-extendable nonplanar bipartite graphs. More precisely, $\mathcal{G}$ consists of all 1-extendable nonplanar bipartite graphs $G$ such that if $G$ contains a 1-extendable nonplanar bipartite graph $H$, then $G$ is isomorphic to $H$. Thus to prove 1.4 we need to show that every member of $\mathcal{G}$ is isomorphic to $K_{3,3}$, a stem, a flower, Uno or Duo. To this end we need two lemmas.

3.4. *Let $G \in \mathcal{G}$, and let $\{e_1, e_2\}$ be a matching of $G$ of size two. Then $G \backslash \{e_1, e_2\}$ is connected.*

*Proof.* The graph $G$ is clearly connected, and hence is 2-connected by 3.1. Suppose for a contradiction that $G \backslash \{e_1, e_2\}$ has two components, say $G_1$ and $G_2$, and let $(A, B)$ be a bipartition of $G$. We shall construct 1-extendable graphs $G'_1, G'_2$ contained in $G$; one of them will contradict the minimality of $G$.

Assume first that $|V(G_1) \cap A| = |V(G_1) \cap B|$. Then by 3.1 each of the sets $V(G_1) \cap A$, $V(G_1) \cap B$, $V(G_2) \cap A$, $V(G_2) \cap B$ contains precisely one vertex incident with $e_1$ or $e_2$. Let $a_1 \in V(G_1) \cap A$, $b_1 \in V(G_1) \cap B$, $a_2 \in V(G_2) \cap A$ and $b_2 \in V(G_2) \cap B$ be those vertices; then one of $e_1, e_2$ has ends $a_1$ and $b_2$, and the other has ends $b_1$ and $a_2$. For $i = 1, 2$ let $G'_i = G_i$ if $a_i$ is adjacent to $b_i$ in $G_i$, and otherwise let $G'_i$ be obtained from $G_i$ by adding an edge joining $a_i$ and $b_i$. It follows from 3.1 that $G'_1$ and $G'_2$ are 1-extendable.

We claim that $G'_1$ and $G'_2$ are both contained in $G$. From the symmetry it suffices to argue for $G'_1$. The claim is obvious if $a_1$ and $b_1$ are adjacent in $G'_1$, and so we may assume that they are not. Let $M_1$ be a perfect matching of $G$ containing one (and hence both) of $e_1$ and $e_2$, and let $M_2$ be a perfect matching of $G$ containing an edge of $E(G) - \{e_1, e_2\}$ incident with $a_2$. (Such an edge exists because $G$ is 2-connected.) The symmetric difference of $M_1$ and $M_2$ contains two paths between $a_1$ and $b_1$; one of these paths, say $P$, includes both $e_1$ and $e_2$, and otherwise is a subgraph of $G_2$. Let $G'''_1$ be obtained from $G_1$ by adding the path $P$; then $G'''_1$ is isomorphic to an even subdivision of $G'_1$ and $G\backslash V(G'''_1)$ has a perfect matching (namely a subset of $M_1$ or $M_2$). This proves that $G'_1$ is contained in $G$, and hence so is $G'_2$. Since $G$ is nonplanar it follows that $G'_1$ or $G'_2$ is nonplanar, contrary to the minimality of $G$. This completes the case when $|V(G_1) \cap A| = |V(G_1) \cap B|$.

Thus we may assume that $|V(G_1) \cap A| > |V(G_1) \cap B|$. Then $|V(G_1) \cap A| = |V(G_1) \cap B| + 1$ by 3.1. It follows that both $e_1$ and $e_2$ have one end in $V(G_1) \cap A$ and the other in $V(G_2) \cap B$. For $i = 1, 2$ let $a_i \in V(G_1) \cap A$ and $b_i \in V(G_2) \cap B$ be the ends of $e_i$, let $G'_1$ be obtained from $G_1$ by adding a new vertex adjacent to $a_1$ and $a_2$, and let $G'_2$ be obtained from $G_2$ by adding a new vertex adjacent to $b_1$ and $b_2$. It follows from 3.1 that $G'_1$ and $G'_2$ are 1-extendable, and an argument similar to the one in the previous paragraph shows that they are both contained in $G$. Since $G$ is nonplanar, one of $G'_1, G'_2$ is nonplanar, contrary to the minimality of $G$. □

3.5. *No member of $\mathcal{G}$ has a magic circuit.*

*Proof.* Let $G \in \mathcal{G}$, let $C$ be a circuit of $G$, and let $e_0 \in E(C)$. Let $M$ be a perfect matching of $G$ containing $e_0$. By 3.3 and 3.4 there exists an edge $e \in E(C) - M$ such that $G\backslash e$ is 1-extendable. It follows that $G\backslash e$ is planar; that is, $e$ is planarizing. Also, since $G\backslash e$ is 1-extendable, both ends of $e$ have degree at least three in $G$. Since $e \neq e_0$ and $e_0$ was arbitrary, we deduce that $C$ is not magic, as desired. □

*Proof of 1.4.* We must show that every member of $\mathcal{G}$ is isomorphic to $K_{3,3}$, a stem, a flower, Uno, or Duo. To this end let $G \in \mathcal{G}$. Then $G$ is clearly connected, and hence it is 2-connected by 3.2. By 3.5 $G$ has no magic circuit, and so by 2.1 it is isomorphic to a subdivision of a graph $H$, where $H$ is one of the graphs listed in 2.1. It follows from the minimality of $G$ that every edge of $H$ is subdivided at most once (that is, is replaced by a path on at most two edges).

Assume first that $H$ is a Möbius ladder, and let $C$ be a circuit of $H$ of length four. Then $C$ contains two rungs and two nonrung edges. Since the circuit of $G$ corresponding to $C$ in $G$ is not magic and every rung of $H$ is nonplanarizing in $H$ we deduce that the nonrung edges of $C$ in $H$ are also edges of $G$ (that is, they are not subdivided in $G$). If some rung of $H$ was not subdivided in $G$, then no rung of $H$ would be subdivided (because $G$ is bipartite), and hence $G = H$. Since $G$ is bipartite, it would have $4n + 2$ vertices for some $n \in \{2, 3 \ldots\}$. By deleting $2(n-1)$ consecutive rungs we see that $G$ contains $K_{3,3}$, contrary to the minimality of $G$. Thus every rung of $H$ is subdivided in $G$, and hence $G$ is a stem.

If $H$ is a biwheel on at least eight vertices, then the same argument applied to the circuits of $H$ of length four not containing the edge joining the two hubs of $H$ shows that $G$ is a flower. Next we assume that $H$ is $K_{3,3}$, $K_{3,3}^+$, $K_{3,3}^{\overset{+}{+}}$, or $K_{3,3}^{++}$. Suppose first that $H$ has an edge $e$ joining two vertices of the same color class (that is, $H \neq K_{3,3}$ and $e$ is one of the "added" edges). We claim that no edge of $K_{3,3}$ adjacent to $e$ in $H$ is subdivided in $G$. Indeed, suppose for a contradiction that $f$ is an edge of $K_{3,3}$ adjacent to $e$ in $H$ such that $f$ is subdivided in $G$. The edges $e$ and $f$ belong to a circuit $C$ of $H$ of length three. Since $G$ is bipartite and $f$ is subdivided in $G$, the circuit of $G$ that corresponds to $C$ is magic, a contradiction. This proves that $f$ is not subdivided in $G$. Since $G$ is bipartite we deduce that $e$ is subdivided in $G$.

Using the facts established above and that the color classes of $G$ have the same size it is easy to see that if $H = K_{3,3}$, then $G$ is isomorphic to $K_{3,3}$ or Bud, if $H = K_{3,3}^+$, then $G$ is isomorphic to Uno, and that if $H = K_{3,3}^{\overset{+}{+}}$, then $G$ is isomorphic to Duo. Similarly, it follows that $H \neq K_{3,3}^{++}$ and $H \neq K_5$. □

## 4. Braces

The objective of this section is to prove 1.5. In a series of lemmas we will show that if a brace contains one of the graphs listed in 1.4, then it contains one of the graphs listed in 1.5. We use the method of "augmenting paths" from network flow theory.

4.1. *Let $H$ be a subgraph of a brace $G$, let $M$ be a perfect matching of $G \backslash V(H)$, let $(A, B)$ be a bipartition of $G$ and let $X$ be a nonempty subset of $A \cap V(H)$ with $|N_H(X)| \leq |X| + 1$ and $N_H(X) \neq B \cap V(H)$. Then there exists an $M$-alternating $H$-path $P$ in $G$ with ends $x \in X$ and $y \in V(H) - N_H(X)$. In particular, $G \backslash V(H \cup P)$ has a perfect matching.*

*Proof.* Let $R$ be the set of all vertices of $G$ that can be reached by an $M$-alternating path with one end in $X$. Thus $X \subseteq R$. If $R \cap B \cap V(H) -$ $N_H(X) \neq \emptyset$ then a corresponding $M$-alternating path is an $H$-path that satisfies the conclusion of the lemma. We may therefore assume that $R \cap B \cap V(H) \subseteq N_H(X)$. But then $|R \cap A - X| = |R \cap B - N_H(X)|$ and so

$$|R \cap A| + 1 = |X| + 1 + |R \cap A - X|$$
$$\geq |N_H(X)| + |R \cap B - N_H(X)| = |N_G(R \cap A)|;$$

hence $N_G(R \cap A) = B$ by 3.1, and thus $N_H(X) = B \cap V(H)$, a contradiction. □

Let $G$ be a graph, let $u$ be a vertex of $G$ of degree two, let $e$ and $f$ be the two edges of $G$ incident with $u$, and let $G'$ be obtained from $G$ by contracting both $e$ and $f$ and deleting parallel edges. We say that $G'$ was obtained from $G$ by a *bicontraction*. We say that a graph $H$ is *weakly contained* in a graph $G$ and that $G$ *weakly contains* $H$ if $G$ has a subgraph $K$ such that $G \backslash V(K)$ has a perfect matching and a graph isomorphic to $H$ can be obtained from $K$ by a sequence of bicontractions. Thus if $H$ is contained in $G$, then it is weakly contained in $G$, but the converse is false. The following is a partial converse.

4.2. *Let $H$ be a graph of maximum degree three. Then $H$ is contained in a graph $G$ if and only if it is weakly contained in $G$.*

*Proof.* The proof is elementary and is omitted. □

We denote the edge of a graph with ends $u$ and $v$ by $uv$ or $u$-$v$. The latter notation will be used when $u$ and $v$ are integers.

Let $G$ be a graph, and let $u, v$ be two vertices of $G$. If $u$ and $v$ are not adjacent we define $G + (u, v)$ to be the graph obtained from $G$ by adding an edge with ends $u, v$; otherwise we define $G + (u, v)$ to be $G$. Now let $G$ be a bipartite graph with bipartition $(A, B)$, let $u \in V(G)$, and let $e$ be an edge of $G$ with ends $a \in A$ and $b \in B$. Let $G'$ be obtained from $G$ by replacing the edge $e$ by a path with vertices $a, b', a', b$ (in the order listed). Thus $G'$ is an even subdivision of $G$. If $u \in A$ we define $G + (u, e)$ to be the graph $G' + (u, b')$, and we say that $b', a'$ (in this order) are the *new* vertices of $G + (u, e)$. If $u \in B$ we define $G + (u, e)$ to be the graph $G' + (u, a')$, and say that $a', b'$ (in this order) are the *new* vertices of $G + (u, e)$. If a graph $H$ is contained in a graph $G$, then $G$ has a central subgraph $K$ isomorphic to an even subdivision of $H$. We say that $K$ is a *model* of $H$ in $G$. Since $K$ is isomorphic to an even subdivision of $H$, there exists a mapping $\phi$ with domain $V(H) \cup E(H)$ such that for all vertices $v, v' \in V(H)$ and all edges $e, e' \in E(H)$,

(i) $\phi(v)$ is a vertex of $K$ and if $v \neq v'$ then $\phi(v) \neq \phi(v')$,
(ii) if $e$ has ends $v$ and $v'$, then $\phi(e)$ is a path in $K$ with ends $\phi(v)$ and $\phi(v')$ and with an even number of vertices,
(iii) if $\phi(v)$ belongs to $\phi(e)$, then it is one of its ends,

(iv) if $e \neq e'$, then $\phi(e)$ and $\phi(e')$ are vertex-disjoint, except possibly for a vertex that is an end of both, and

(v) $V(K) = \bigcup_{u \in V(H)} \{\phi(u)\} \cup \bigcup_{f \in E(H)} V(\phi(f))$ and $E(K) = \bigcup_{f \in E(H)} E(\phi(f))$.

We say that $\phi$ is an *embedding* of $H$ into $G$, and we write $\phi(H) = K$.

Let $H$ be a bipartite graph with bipartition $(A, B)$, and let $X \subseteq A$. If $|N_H(X)| \leq |X| + 1$, $N_H(X) \neq B$, and $H$ is contained in a brace $G$, then the analogous inequality does not hold in $G$ by 3.1, and thus $G$ contains an "augmentation" of $H$. Let us make this precise. We define $\Phi_H(X)$ to be the set of all edges of $H$ with ends $u$ and $v$, where $u \in N_H(X)$ and $v \notin X$. Let $H'$ be a subdivision of $H$. For $v \in V(H)$ we define $\langle v \rangle_X^{H'} = v$. If $e \in \Phi_H(X)$ has ends $u$ and $v$, where $u \in N_H(X)$ and $v \notin X$, then the edge $e$ corresponds to a path $P$ in $H'$ with ends $u$ and $v$. We define $\langle e \rangle_X^{H'}$ to be the edge of $P$ incident with $v$.

Now let $k \geq 1$ be an integer, let $e_1, e_2, \ldots, e_{k-1} \in \Phi_H(X)$, and let $e_k \in \Phi_H(X)$, or let $e_k$ be a vertex of $B - N_H(X)$, or an edge of $H$ not incident with a vertex of $N_H(X)$. Let $a_0 \in X$. We define $H_0 = H$, and for $i = 1, 2, \ldots, k$ we define recursively $H_i = H_{i-1} + (a_{i-1}, \langle e_i \rangle_X^{H_{i-1}})$, and if $e_i \notin V(H)$ let $b_i, a_i$ be the new vertices of $H_i$. We say that $H_k$ is a *partial $X$-augmentation* of $H$. If $e_k \notin \Phi_H(X)$ we say that $H_k$ is an *$X$-augmentation* of $H$. In either case we say that $H_k$ is *determined* by $a_0, e_1, e_2, \ldots, e_k$. See Figure 4 below, where $a_0 = w$. Later on we will also consider $H$ to be a partial $X$-augmentation of itself determined by $a_0$.

Now let $G$ be a graph, and let $\phi$ be an embedding of $H$ into $G$. We say that an embedding $\phi_k$ of $H_k$ into $G$ *extends* $\phi$ if $\phi(v) = \phi_k(v)$ for every $v \in V(H) - X$ and $\phi(H) = \phi_k(H')$, where $H'$ is obtained from $H_k$ by deleting the edges $a_{i-1}b_i$ for all $i = 1, 2, \ldots, k$. Let $M$ be a perfect matching of $G \backslash V(\phi(H))$. We say that $\phi_k$ is *$M$-compatible* if for all $i = 1, 2, \ldots, k$, the path $\phi_k(a_{i-1}b_i)$ is $M$-alternating. The following result is our main tool in the proof of 1.5.

4.3. *Let $H$ be a bipartite graph with bipartition $(A, B)$, let $H$ be contained in a brace $G$, and let $X \subseteq A$ be a nonempty set of vertices with $|N_H(X)| \leq |X| + 1$ and $N_H(X) \neq B$. Then some $X$-augmentation of $H$ is weakly contained in $G$. Moreover, if the degree in $H$ of every vertex of $X$ is at most two, then some $X$-augmentation is contained in $G$.*

*Proof.* We only prove the second statement; the proof of the first is similar. Let $\phi$ be an embedding of $H$ into $G$, let $K = \phi(H)$, let $M$ be a perfect matching of $G \backslash V(K)$, and let $(A', B')$ be the bipartition of $G$ with $\phi(A) \subseteq A'$. Let $X_1'$ be the set of all vertices of $K$ that belong to $V(\phi(e)) \cap A'$ for some edge $e$ of $H$ incident with a vertex of $X$. Since $|N_H(X)| \leq |X| + 1$ we deduce that $|N_K(X_1')| \leq |X_1'| + 1$.

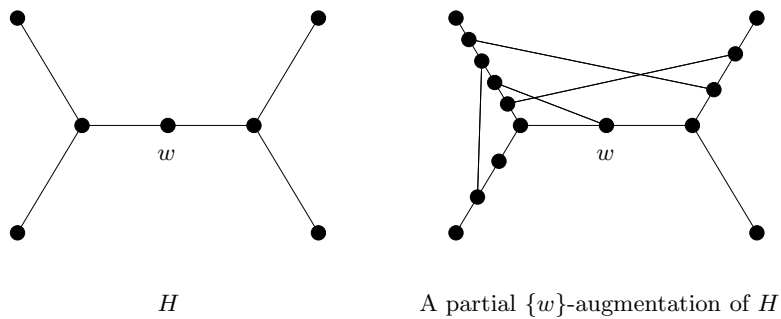

$H$          A partial $\{w\}$-augmentation of $H$

Figure 4

Let $X'_2$ be the set of all vertices of $K$ that can be written as $\phi_k(a_k)$, where $k \geq 1$ is an integer, $\phi_k$ is an $M$-compatible embedding of a partial $X$-augmentation $H_k$ of $H$ into $G$ extending $\phi$, $H_k$ is not an $X$-augmentation, $H_k$ is determined by $a_0, e_1, e_2, \ldots, e_k$ and the vertices $a_1, a_2, \ldots, a_k$ are as in the definition of partial $X$-augmentation. Since $H_k$ is not an $X$-augmentation it follows that $X'_2 \subseteq \bigcup_{e \in \Phi_H(X)} V(\phi(e)) \cap A'$. It is easy to see that if $e \in \Phi_H(X)$ has ends $u$ and $v$, where $u \in N_H(X)$ and $v \notin X$, and if $a \in V(\phi(e)) \cap X'_2$, then every vertex of the subpath of $\phi(e)$ between $a$ and $\phi(u)$ that belongs to $A'$ belongs to $X'_2$. We deduce that $|N_K(X'_2) - N_K(X'_1)| \leq |X'_2|$. Now let $X' = X'_1 \cup X'_2$; then

$$|N_K(X')| \leq |N_K(X'_1)| + |N_K(X'_2) - N_K(X'_1)| \leq |X'_1| + 1 + |X'_2| = |X'| + 1.$$

Moreover, $N_K(X') \neq B'$, because $N_H(X) \neq B$.

By 4.1 there exists an $M$-alternating $K$-path $P$ in $G$ with ends $x \in X'$ and $y \in V(K) - N_K(X')$. If $y \in \phi(V(H))$, let $z \in B - N_H(X)$ be such that $\phi(z) = y$; otherwise let $z$ be an edge of $H$ such that $y \in V(\phi(z))$. Since $x \in X'$ there exist $k \geq 0$, $H_k$, $a_0, e_1, e_2, \ldots, e_k, a_k$ and $\phi_k$ such that $x = \phi_k(a_k)$, $H_k$ is a partial $X$-augmentation of $H$ determined by $a_0, e_1, e_2, \ldots, e_k$, the vertices $a_1, a_2, \ldots, a_k$ are as in the definition of $X$-augmentation, and $\phi_k$ is an $M$-compatible embedding of $H_k$ into $G$ extending $\phi$. Let $t$ be the minimum integer such that either $t \geq k$, or $t \in \{0, 1, \ldots, k-1\}$ and $V(P) \cap V(\phi_k(a_t b_{t+1})) \neq \emptyset$. Since $\phi_k$ is $M$-compatible and $P$ is $M$-alternating, it follows that $\phi_k(H_k) \cup P$ contains the partial $X$-augmentation $H'$ of $H$ determined by $a_0, e_1, e_2, \ldots, e_t, z$ if $k \geq 1$ and by $x, z$ if $k = 0$. Moreover, there exists an $M$-compatible embedding of $H'$ into $G$ extending $\phi$. (To see this, notice that if $t = k$, then $\phi_k(H_k) \cup P$ is as desired. Otherwise delete $\phi_k(a_j b_{j+1})$ for all $j = t+1, t+2, \ldots, k-1$, and suitable subpaths of $P$ and $\phi_k(a_t b_{t+1})$.) We claim that $z \notin \Phi_H(X)$. Indeed, suppose for a contradiction that $z$ is an edge that belongs to $\Phi_H(X)$, and let $y'$ be the neighbor of $y$ on the path $\phi(z)$ chosen so that $\phi_k(u), y', y, \phi_k(v)$ occur on $\phi(z)$ in the order listed, where $u \in N(X)$ and $v \notin X$ are the ends of $z$. Then $y' \in X_2'$ by the definition of $X_2'$, contrary to the fact that $y \notin N_K(X')$. Thus $z \notin \Phi_H(X)$, and hence $H'$ is an $X$-augmentation, as desired. □

4.4. *Let $G$ be a brace, and let $H$ be a bipartite graph contained in $G$. Let $(A,B)$ be a bipartition of $H$ with $|B| \geq 3$, and let $w \in A$ be a vertex of $H$ of degree two. Then either*

(i) *there exists a vertex $v \in B$ not adjacent to $w$ in $H$ such that $G$ weakly contains $H + (w,v)$, or*

(ii) *there exists an edge $e$ of $H$ incident with a neighbor of $w$ but not with $w$ itself such that $G$ contains $H + (w,e)$.*

*Proof.* This follows immediately from 4.3 applied to $X = \{w\}$, because $|B| \geq 3$. □

4.5. *Let $G$ be a 2-extendable bipartite graph containing Uno. Then $G$ contains either $K_{3,3}$ or Bud.*

*Proof.* Let $H$ be Uno with its vertices numbered as in Figure 2. By 4.2 it suffices to show that $G$ weakly contains $K_{3,3}$ or Bud. By 4.4 and taking symmetry into account $G$ weakly contains either $H + (10,7)$ or $H + (10, 2\text{-}4)$. But $H + (10,7)$ contains $K_{3,3}$ (delete 2-6 and 3-5), and $H + (10, 2\text{-}4)$ contains Bud (delete 2-10 and 3-4). Thus $G$ weakly contains $K_{3,3}$ or Bud, as desired. □

4.6. *Let $G$ be a 2-extendable bipartite graph containing Duo. Then $G$ contains $K_{3,3}$ or Bud.*

*Proof.* Let $H$ be Duo with its vertices numbered as in Figure 2. By 4.4 and taking symmetry into account $G$ weakly contains $H + (7,8)$ or $H + (7,1)$, or it contains $H + (7, 2\text{-}5)$ or $H + (7, 2\text{-}4)$. The first graph contains $K_{3,3}$ (delete the vertices 1 and 4), the second contains $K_{3,3}$ (delete the vertices 6 and 8), the third contains Bud (delete 2-7, 3-5 and 1-6), and the fourth contains Uno (delete 2-7 and 3-4). Thus either $G$ weakly contains $K_{3,3}$ (in which case it contains $K_{3,3}$ by 4.2), or $G$ contains Bud or Uno, in which case it contains $K_{3,3}$ or Bud by 4.5. □

We need two refinements of 4.3, and we now introduce the first of the two. Let $H_k$ be an $X$-augmentation of $H$ determined by $a_0, e_1, e_2, \ldots, e_k$. We say that $H_k$ is *weakly reduced* if the following two conditions hold for all $i = 1, 2, \ldots, k-2$.

(R1) The edges $e_i$ and $e_{i+1}$ are distinct, and

(R2) no vertex of $N_H(X)$ of degree three in $H$ is an end of both $e_i$ and $e_{i+1}$.

4.7. *Let $H$ be a bipartite graph with bipartition $(A, B)$, let $H$ be contained in a brace $G$, and let $X \subseteq A$ be a nonempty set of vertices of $H$ of degree two such that $|N_H(X)| \leq |X| + 1$ and $N_H(X) \neq B$. Then some weakly reduced $X$-augmentation of $H$ is contained in $G$.*

*Proof.* By 4.3 we may choose the minimum integer $k \geq 1$ and elements $a_0, e_1, e_2, \ldots, e_k$ such that $G$ contains an $X$-augmentation $H_k$ of $H$ determined by $a_0, e_1, e_2, \ldots, e_k$. We claim that $H_k$ is as desired. Indeed, suppose first that $e_i = e_{i+1}$ for some $i \in \{1, 2, \ldots, k-2\}$. Let $a_i, b_i, a_{i+1}, b_{i+1}$ be as in the definition of $X$-augmentation, and let $H'$ be obtained from $H_k$ by deleting the edges $a_i b_i$ and $a_{i+1} b_{i+1}$. Then $H'$ is isomorphic to the subdivision of the $X$-augmentation of $H$ determined by $a_0, e_1, \ldots, e_{i-1}, e_{i+2}, \ldots, e_k$, contrary to the choice of $k$, because $H'$ is contained in $G$.

Suppose now that $v \in N(X)$ has degree three, and that it is an end of both $e_i$ and $e_{i+1}$ for some $i \in \{1, 2, \ldots, k-2\}$. Let $w$ be the unique neighbor of $v$ in $H$ that belongs to $X$. For $j = 1, 2, \ldots, k$ let $a_j, b_j$ be as in the definition of $X$-augmentation. Let $\phi_k$ be an embedding of $H_k$ into $G$, let $\phi_k(H_k) = K$, and let $P$ be the path of $K$ corresponding to the edge $e_i$ of $H$. Let $K'$ be obtained from $K$ by deleting the edges and interior vertices of the subpath of $P$ with ends $\phi_k(v)$ and $\phi_k(a_i)$, and for $j = 0, 1, \ldots, i-1$ deleting the edges and interior vertices of the path $\phi_k(a_j b_{j+1})$. Since $w$ has degree two in $H$ we deduce that $K'$ is isomorphic to an even subdivision of the $X$-augmentation of $H$ determined by $w, e_{i+2}, e_{i+3}, \ldots, e_k$, contrary to the choice of $k$. □

4.8. *Let $G$ be a brace containing a flower. Then $G$ contains Bud.*

*Proof.* For $k = 2, 3, \ldots$ let $F_k$ be the flower on $4k + 2$ vertices. Thus $F_2$ is Bud. We will prove that if $G$ contains $F_k$ for some integer $k \geq 3$, then it contains $F_j$ for some $j = 2, 3, \ldots, k-1$. Let the vertices of $F_k$ be $u_0, u_1, u_2, \ldots, u_{2k}, u_{2k+1}, v_1, v_2, \ldots, v_{2k}$ in such a way that the vertices $u_1, u_2, \ldots, u_{2k}$ form a circuit (in the order listed), and for $i = 1, 2, \ldots, k$, the neighbors of $v_{2i}$ are $u_{2i}$ and $u_{2k+1}$, and the neighbors of $v_{2i-1}$ are $u_{2i-1}$ and $u_0$. Let us say that a graph is *good* if it contains $F_j$ for some $j = 2, 3, \ldots, k-1$.

(1) *Let $i \in \{1, 2, \ldots, k\}$, and let $x$ be $v_{2i}$, or an edge of $F_k$ incident with $v_{2i}$. Then $F_k + (v_1, x)$ is good.*

To prove (1) we notice that $F_k + (v_1, x)$ contains $F_k + (v_1, v_{2i})$, and so we may assume that $x = v_{2i}$. We first assume that $k = 3$ and $i = 2$. Then $F_k + (v_1, x)$ contains $F_2$ (delete $u_5$, $v_5$, $u_6$ and $v_6$). Thus from the symmetry we may assume that $i \geq 3$. But then $F_k + (v_1, v_{2i})$ contains $F_{i-1}$ (delete $u_j$ and $v_j$ for $j = 2i + 1, \ldots, 2k$, and the edges $u_0 v_1$, $v_{2i} u_{2k+1}$ and $u_{2k} u_1$). This proves (1).

(2) Let $i \in \{2, \ldots, k\}$, and let $x$ be $u_{2i-1}$, or an edge incident with $u_{2i-1}$, or an edge incident with $v_{2i-1}$. Then $F_k + (v_1, x)$ is good.

To prove (2) we may assume from the symmetry that $i \geq 3$. Then $F_k + (v_1, x)$ contains $F_{i-1}$ (delete $u_j$ and $v_j$ for $j = 2i, 2i+1, \ldots, 2k$ and the edge $u_0 v_1$), and (2) follows.

Let $H = F_k + (v_1, u_1 u_2)$, and let $a, b$ be the new vertices of $H$. Thus $b$ has degree two in $H$. Let $H'$ be obtained from $H$ by deleting the edge $v_1 a$. From (2) we deduce the following.

(3) Let $i \in \{2, 3, \ldots, k\}$ and let $x$ be $v_{2i}$ or an edge incident with $v_{2i}$. Then $H' + (b, x)$ is good.

To prove (3) we first notice that $H' + (b, x)$ contains $H' + (b, v_{2i})$. But $H' + (b, v_{2i})$ is isomorphic to $F_k + (v_1, u_{2i-1} u_{2i})$, which is good by (2). This proves (3).

(4) Let $i = 3, 4, \ldots, k$, and let $x$ be $u_{2i-1}$, or an edge incident with $u_{2i-1}$, or an edge incident with $v_{2i-1}$. Then $H' + (b, x)$ is good.

Indeed, $H' + (b, x)$ contains $F_{i-1}$ (delete $u_j$ and $v_j$ for $j = 1$ and $j = 2i, 2i+1, \ldots, 2k$). This proves (4).

We are now ready to prove that if a brace $G$ contains $F_k$ for some integer $k \geq 3$, then it contains $F_j$ for some $j = 2, 3, \ldots, k-1$. To this end let $k \geq 3$ be an integer, and assume that a brace $G$ contains $F_k$. By 4.7 a weakly reduced $\{v_1\}$-augmentation $K$ of $F_k$ is contained in $G$. Let $K$ be determined by $v_1, e_1, e_2, \ldots, e_t$, where $t \geq 1$. Let us first assume that $t = 1$; then either $e_1 \in \{u_3, u_5, \ldots, u_{2k-1}, v_2, v_4, \ldots, v_{2k}\}$, or $e_1$ is an edge of $F_k$ incident with one of those vertices. By (1) and (2) $K$ is good, and hence so is $G$, as desired. Thus we may assume that $t > 1$, and hence $e_1$ is an edge of $F_k$ incident with $u_0$ or $u_1$, but not with $v_1$. If $e_1$ is the edge $u_0 u_{2k+1}$, then $G$ contains Bud, because $F_k + (v_1, u_0 u_{2k+1})$ does (delete $u_i$ and $v_i$ for $i = 4, 5, \ldots, 2k-2$ and the edge $u_0 v_1$). By (2) we may assume that $e_1$ is not incident with $u_0$, and so from the symmetry we may assume that $e_1$ is the edge $u_1 u_2$ of $F_k$. If $e_2$ is the edge $u_0 u_{2k+1}$, then $G$ contains $K_{3,3}$, because $H + (b, u_0 u_{2k+1})$ does (delete $u_i$ and $v_i$ for $i = 3, 4, \ldots, 2k-1$). By (R1) and (R2) $e_2$ is not an edge of $F_k$ incident with $u_1$, and so by (3) and (4) we may assume that $e_2$ is $v_2$, or $u_3$, or an edge incident with $v_2$ or $u_3$. In each case $K$ contains $F_{k-1}$. It suffices to check this for $e_2 \in \{v_2, u_3\}$. If $e_2 = v_2$ delete $u_0 v_1$, $v_2 u_{2k+1}$, $a u_2$ and $u_1 b$; if $e_2 = u_3$ delete $u_2, v_2$ and the edges $u_0 v_1$ and $u_1 b$. □

We now introduce the second refinement of 4.3. Let $H$ be a bipartite graph with bipartition $(A, B)$, let $X \subseteq A$, let $k \geq 1$ be an integer, and let $H_k$ be an $X$-augmentation of $H$ determined by $a_0, e_1, e_2, \ldots, e_k$. We say that $H_k$ is *reduced* if $H_k$ is weakly reduced, and for all $i = 1, 2, \ldots, k-1$, letting $u_i$ be the unique end of $e_i$ that belongs to $N_H(X)$, the following conditions hold.

(R3) The vertex $u_i$ has degree at least three in $H$;

(R4) if $u_1$ has degree three in $H$ and is adjacent to $a_0$ in $H$, then $a_0$ is the only neighbor of $u_1$ in $H$ that belongs to $X$;

(R5) $e_k \in B - N(X)$;

(R6) if $k = 2$, if $a_0$ is adjacent to $u_1$ in $H$ and if $e_2$ is a vertex of $H$ adjacent to a neighbor $v$ of $u_1$ in $H$, then either $u_1$ has degree at least four in $H$, or $v$ has degree at least four in $H$, or $v$ has degree three in $H$ and the neighbor of $v$ in $H$ not in $\{u_1, e_2\}$ belongs to $N_H(X)$; and

(R7) if $k \geq 3$, if $a_0$ has degree two in $H$ and $u_1$ and $u_2$ are its neighbors in $H$, and if $u_2$ has degree three in $H$, then $a_0$ is the only neighbor of $u_2$ in $H$ that belongs to $X$.

4.9. *Let $H$ be a bipartite graph with bipartition $(A, B)$, let $H$ be contained in a brace $G$, and let $X \subseteq A$ be a nonempty set such that $|N_H(X)| \leq |X| + 1$ and $N_H(X) \neq B$. Then some reduced $X$-augmentation of $H$ is weakly contained in $G$.*

*Proof.* By 4.3 we may choose the minimum integer $k \geq 1$ and elements $a_0, e_1, e_2, \ldots, e_k$ such that $G$ weakly contains the $X$-augmentation $H_k$ of $H$ determined by $a_0, e_1, e_2, \ldots, e_k$, and subject to that, $|V(H_k)|$ is minimum. It follows that (R5) holds. The argument used in the proof of 4.7 shows that $H_k$ satisfies (R1) and (R2). For $i = 1, 2, \ldots, k-1$ let $u_i$ be the end of $e_i$ that belongs to $N_H(X)$. To prove (R3) let $i \in \{1, 2, \ldots, k-1\}$. If $u_i$ does not have degree at least three, then it has a unique neighbor in $H$, say $x$, that belongs to $X$. It follows that the $X$-augmentation of $H$ determined by $x, e_{i+1}, e_{i+2}, \ldots, e_k$ is weakly contained in $G$, contrary to the choice of $k$. This proves that $H_k$ satisfies (R3).

To prove that $H_k$ satisfies (R4) suppose that $u_1$ has degree three, and that $a_0, x \in X$ are two distinct neighbors of $u_1$. Let $H'$ be obtained from $H_k$ by deleting the edge $u_1 a_0$; then $H'$ is weakly contained in $G$, and it weakly contains the $X$-augmentation of $H$ determined by $x, e_2, e_3, \ldots, e_k$, contrary to the choice of $k$. This proves that $H_k$ satisfies (R4).

To prove that $H_k$ satisfies (R6) let $k = 2$, let $a_0$ be adjacent to $u_1$ in $H$, and let $e_2$ be a vertex of $H$ adjacent to a neighbor $v$ of $u_1$. By (R3) we may assume that $u_1$ has degree three in $H$. Let $b_1, a_1$ be the new vertices of $H_1 = H + (a_0, e_1)$. From the symmetry between $u_1$ and $b_1$ we may assume that $v$ is an end of $e_1$. If $v$ has degree two in $H$, then $H_2$ weakly contains the $X$-augmentation of $H$ determined by $a_0, e_2$, contrary to the choice of $k$. We may therefore assume that $v$ has degree three in $H$, and that $y \notin N_H(X)$, where $y$ is the neighbor of $v$ in $H$ not in $\{u_1, e_2\}$. By deleting the edge $ve_2$ of $H_2$ we see that $H_2$ weakly contains the $X$-augmentation of $H$ determined by $a_0, y$, contrary to the choice of $k$. This proves that $H_k$ satisfies (R6).

It remains to prove that $H_k$ satisfies (R7). Let $a_0$ have degree two in $H$, let $u_1$ and $u_2$ be its neighbors in $H$, let $u_2$ have degree three in $H$, and suppose that $u_2$ has a neighbor $x \in X - \{a_0\}$ in $H$. Then $x$ is not incident with $e_2$, because $e_2 \in \Phi_H(X)$. Let $b_1, a_1$ be the new vertices of $H_1 = H + (a_0, e_1)$. Let $H'$ be obtained from $H_k$ by deleting the edges $a_1 b_1$ and $a_0 u_2$; then $H'$ is weakly contained in $G$ and weakly contains the $X$-augmentation of $H$ determined by $x, e_3, e_4, \ldots, e_k$, contrary to the choice of $k$. This proves that $H_k$ satisfies (R7), and hence completes the proof of the lemma. □

For convenience we state the following corollary of 4.9. Let $H$ be a bipartite graph with bipartition $(A, B)$, let $w \in A$ have degree two, let $u_1, u_2$ be the two neighbors of $w$. If $v \in B - \{u_1, u_2\}$ we say that $H + (w, v)$ is a *w-extension of $H$ of the first kind*. Let $e \in E(H)$ be an edge of $H$ incident with $u_1$ but not with $w$, let $u_1$ have degree at least three, and let $b, a$ be the new vertices of $H' = H + (w, e)$. Let $v$ be a vertex in $B - \{u_1, u_2\}$ with the property that if $v$ is adjacent in $H$ to a neighbor $v'$ of $u_1$, then $u_1$ has degree at least four, or $v'$ has degree at least four, or $v'$ has degree three and is adjacent to $u_2$. We say that $H' + (a, v)$ is a *w-extension of $H$ of the second kind*. Thirdly, let $f$ be an edge incident with $u_2$ but not with $w$, let $u_1$ and $u_2$ have degree at least three, and let $H'$ be as above. We say that $H' + (a, f)$ is a *w-extension of $H$ of the third kind*. Finally, we say that a graph is a *w-extension* of $H$ if it is a $w$-extension of the first, second or third kind.

4.10. *Let $H$ be a bipartite graph with bipartition $(A, B)$ contained in a brace $G$, and let $w \in A$ be a vertex of $H$ of degree two. If both neighbors of $w$ in $H$ have degree at most three and $|B| \geq 3$, then some $w$-extension of $H$ is weakly contained in $G$.*

*Proof.* By 4.9 $G$ weakly contains a reduced $\{w\}$-augmentation $H'$ of $H$. Let $H'$ be determined by $w, e_1, e_2, \ldots, e_k$. If $k = 1$, then $H'$ is a $w$-extension of $H$ of the first kind by (R5). If $k = 2$, then $H'$ is a $w$-extension of $H$ of the second kind by (R3), (R5) and (R6). Finally, if $k \geq 3$, then $H'$ contains a $w$-extension of $H$ of the third kind by (R1), (R2) and (R3). □

4.11. *Let $G$ be a brace containing a stem. Then $G$ contains $K_{3,3}$ or Bud.*

*Proof.* For $k = 2, 3, \ldots$, let $S_k$ denote the stem on $6k$ vertices, and let $S_1$ denote Bud. For $k \geq 2$ let the vertices of $S_k$ be $u_1, u_2, \ldots, u_{4k}, v_1, v_2, \ldots, v_{2k}$ in such a way that $u_1, u_2, \ldots, u_{4k}$ form a circuit in the order listed, and for $i = 1, 2, \ldots, 2k$, the neighbors of $v_i$ are $u_i$ and $u_{2k+i}$. Let us say that a graph is *good* if it contains $K_{3,3}$ or $S_j$ for some $j = 1, 2, \ldots, k-1$.

(1) If $i \in \{1, 2, \ldots, k\}$, then $S_k + (v_1, v_{2i})$ is good. If $i \in \{2, 3, \ldots, 2k\} - \{k+1\}$, then $S_k + (v_1, u_{2i-1})$ is good.

To prove (1) we first notice that $S_k + (v_1, v_2)$ contains $S_{k-1}$ (if $k = 2$ delete $u_1v_1$ and $v_2u_6$; otherwise delete $u_1v_1$, $u_2v_2$ and $u_{2k+1}u_{2k+2}$). So let $x = v_{2i}$, where $i \in \{2, \ldots, k-1\}$, or $x = u_{2i-1}$, where $i \in \{2, 3, \ldots, 2k\} - \{k+1\}$. Then the graph $S_k + (v_1, x)$ contains $K_{3,3}$ (use the edges $v_1u_1$, $v_1u_{2k+1}$, $u_1u_2$, $u_2v_2$, $v_2u_{2k+2}$, $u_{2k+1}u_{2k+2}$, $u_{2k}u_{2k+1}$, $v_{2k}u_{2k}$, $v_{2k}u_{4k}$, $u_1u_{4k}$, $v_1x$, and a suitable path between $u_{2k}$ and $u_{2k+2}$ that contains $x$). This proves (1).

(2) If $x = u_{2i}$ for some $i \in \{2, 3, \ldots, 2k-1\} - \{k, k+1\}$, or $x = v_{2i-1}$ for some $i = 2, 3, \ldots, k$, then $S_k + (u_1, x)$ is good.

To prove (2) we may assume, by (1), that $x = u_{2i}$ for some $i \in \{2, 3, \ldots, 2k-1\} - \{k, k+1\}$. Moreover, from the symmetry we may assume that $i < k$. By deleting $v_j$ and $u_j$ for $j = 2, 3, \ldots, 2i - 1$ we see that $S_k + (u_1, x)$ contains $S_{k-i+1}$, as desired. This proves (2).

We are now ready to complete the proof. To this end let $k \geq 2$ be an integer, and let $G$ contain $S_k$; we will show that $G$ is good. By 4.10 $G$ weakly contains a $\{v_1\}$-extension $H$ of $S_k$. If $H$ is of the first kind, then the lemma holds by (1). Next assume that $H$ is of the second kind. From the symmetry we may assume that $H = S_k + (v_1, u_1u_2) + (a, x)$, where $b, a$ are the new vertices of $S_k + (v_1, u_1u_2)$, and $x \in \{u_5, u_7, \ldots, u_{2k-1}, v_4, v_6, \ldots, v_{2k-2}\}$. (Notice that because of the symmetry between $u_1$ and $b$ we may assume that $x \notin \{u_{2k+3}, u_{2k+5}, \ldots, u_{4k-3}\}$.) Then $H$ weakly contains $S_k + (u_2, x)$, but the latter graph is good by (2) because of symmetry. We may therefore assume that $H$ is of the third kind. From the symmetry between $v_1u_1$ and $v_1b$ we may assume that $H = S_k + (v_1, u_1u_2) + (a, u_{2k+1}u_{2k+2})$. But then $H$ contains Bud (delete the vertices $u_j$ and $v_j$ for $j = 3, 4, \ldots, 2k-1$ and the edge $u_1v_1$), and hence so does $G$, as desired. □

Let us introduce the following convention. If $G$ is a bipartite graph with $V(G) = \{1, 2, \ldots, n\}$, $v \in V(G)$ and $e \in E(G)$, then the new vertices of $G + (v, e)$ are $n + 1$ and $n + 2$. Thus $n + 2$ has degree two in $G + (v, e)$. We define Superbud as Bud$+(10, 5\text{-}6) + (9, 12)$. See Figure 2.

4.12. *Let $G$ be a brace containing Bud. Then $G$ contains $K_{3,3}$ or Superbud.*

*Proof.* Let $G$ be a brace, let $H$ be Bud with vertices numbered as in Figure 2, and let $X = \{7, 9\}$. By 4.9 $G$ weakly contains a reduced $X$-augmentation $H'$ of $H$. Let $H'$ be determined by $a_0, e_1, e_2, \ldots, e_k$, where $k \geq 1$ is an integer. From the symmetry we may assume that $a_0 = 7$. By (R5) $e_k \in \{8, 10\}$. If $k = 1$, then $H'$ contains $K_{3,3}$ (delete the edge 1-2), and so we may assume that $k > 1$. Since $e_1$ is not 1-1 by (R4), we may assume from the symmetry that $e_1$ is either 3-4, or 4-5.

Assume first that $e_1$ is the edge 3-4. If $k = 2$, then $H'$ is isomorphic to Superbud, and so we may assume in this case that $k > 2$. By (R1) $e_2$ is not 3-4, by (R2) $e_2$ is not 3-6, by (R7) $e_2$ is not 1-2, and so we may assume from the symmetry of $H + (7, 3\text{-}4)$ that $e_2$ is 4-5. (The symmetry fixes $1, 2, 5, 7, 9, 12$ and exchanges the pairs $(3, 11)$, $(4, 6)$, and $(8, 10)$). But $H + (7, 3\text{-}4) + (12, 4\text{-}5)$ contains $K_{3,3}$ (delete 1-2 and 4-13), which completes the analysis in the case when $e_1$ is 3-4.

We may therefore assume that $e_1$ is the edge 4-5. Then $e_2 \neq 5\text{-}6$ by (R2). If $k = 2$ and $e_2 = 8$, then $H'$ is isomorphic to Superbud; if $e_2 \in \{3\text{-}4, 3\text{-}6\}$, or $k = 2$ and $e_2 = 10$, then $H'$ contains $K_{3,3}$. To see this when $e_2 = 10$ delete 2-10 and 5-6; the case $e_2 = 3\text{-}6$ is similar, and if $e_2 = 3\text{-}4$ then delete the vertices 6 and 10. We may therefore assume that $k > 2$ and that $e_2$ is the edge 1-2. If $k = 3$, then either $e_3 = 8$, in which case $H'$ contains $K_{3,3}$ (delete 2-8, 5-6 and 13-14), or $e_3 = 10$, in which case $H'$ contains Superbud (delete 2-10). Thus we may assume that $k > 3$. If $e_3$ is 3-4 or 3-6, then $H'$ weakly contains $L = H + (7, 4\text{-}5) + (12, 1\text{-}2) + (14, 3)$, but the graph $L$ contains $K_{3,3}$ (delete 3-6, 5-12 and 2-13). If $e_3$ is 4-5 then $H'$ contains $K_{3,3}$ (delete the vertices 6 and 10, and the "new" edge incident with 7). Thus we may assume that $e_3$ is 5-6. If $e_4 = 10$, then $H'$ weakly contains the graph $L$ (delete 6-10), and hence $H'$ contains $K_{3,3}$. If $e_4 = 8$ or $e_4$ is 4-5, then $H'$ contains $K_{3,3}$. To see this delete 3-4, 7-11, 6-15 and 5-16. If $e_4$ is 3-6 or 3-4, then $H'$ contains Superbud (delete the vertices 6 and 10). Thus we may assume that $e_4$ is 1-2, but then $H'$ contains $K_{3,3}$ (delete vertices 11 and 12, and the edge 3-6).  □

4.13. *Let $G$ be a brace containing Superbud. Then $G$ contains $K_{3,3}$, the Heawood graph or Rotunda.*

*Proof.* Let $H$ be Superbud with vertices numbered as in Figure 2. By 4.2 it suffices to show that $G$ weakly contains $K_{3,3}$, the Heawood graph, or Rotunda. By 4.10 and the symmetry of $H$, $G$ weakly contains one of $H + (7, 8)$, $H + (7, 5)$, $H + (7, 12)$, or $H + (7, 3\text{-}4) + (14, 1)$. The first two graphs contain $K_{3,3}$ (in the first case delete 1-7, 2-8 and 3-4, and in the second case delete the vertices 4 and 8 and the edge 1-7), and so does the last (delete 5, 9, 11, 12). Thus we may assume that $G$ weakly contains $H+(7,12)$; that is, $G$ contains $H_1 = H+(7,12)$, $H_2 = H+(7,12\text{-}6)$, $H_3 = H+(7,12\text{-}9)$, or $H_4 = H+(7,12\text{-}11)$. By 4.10 there is an integer $i \in \{1,2,3,4\}$ such that $G$ weakly contains an 8-extension $K$ of $H_i$. If $K$ is of the first kind, then by the previous containments and symmetry between 7 and 8 we may assume that $G$ weakly contains $H+(7,12)+(8,11)$ or $H+(7,12\text{-}11)+(8,14)$. The first graph contains $K_{3,3}$ (delete 1-9, 2-8, 5-11 and 11-12), and the second is isomorphic to the Heawood graph. This completes the case when $K$ is of the first kind.

Assume now that $K$ is of the second or third kind. Then either $K$ weakly contains $H+(8,3\text{-}4)+(14,2)$, or $i > 1$ and $K = H_i + (8,3\text{-}4) + (16,14)$. Since $H+(8,3\text{-}4)+(14,2)$ is isomorphic to $H+(7,3\text{-}4)+(14,1)$, and hence contains $K_{3,3}$, we may assume that $i > 1$ and that $K = H_i + (8,3\text{-}4) + (16,14)$. If $i = 2$, then $K$ contains $K_{3,3}$ (delete 1-2, 5-11, 6-10, 7-13 and 9-12), if $i = 3$ then $K$ is isomorphic to Rotunda, and if $i = 4$ then $K$ contains $K_{3,3}$ (delete the vertices 1 and 7, and edges 8-15 and 12-14). □

Theorem 1.5 now follows from 1.4, 4.5, 4.6, 4.8, 4.11, 4.12 and 4.13.

## 5. Rotunda

Let $G$ be a brace not containing $K_{3,3}$, let $K$ be a model of Rotunda in $G$, and let $S$ be the set of four vertices of $K$ corresponding to the center of Rotunda. The objective of this section is to show that in those circumstances each component of $K\backslash S$ belongs to a different component of $G\backslash S$. We say that a path $P$ in $G$ is a $(K,S)$-*jump* if $P$ is a $K$-path in $G$ such that its ends belong to different components of $K\backslash S$ (and hence $P$ is disjoint from $S$), and such that $G\backslash V(K \cup P)$ has a perfect matching. We shall need to show that there is no $(K,S)$-jump in $G$, but before that we need three lemmas, which require some definitions.

Let $H$ be a bipartite graph with bipartition $(A,B)$, let $u \in A$ be a vertex of $H$ of degree three, and let $u_1, u_2, u_3$ be its neighbors. Let $H'$ be obtained from $H$ by replacing, for $i = 1,2$, the edge $uu_i$ by a path with vertex-set $u, b_i, a_i, u_i$ (in order), and let $H'' = H' + (a_1, b_2) + (b_1, a_2)$. We say that $H''$ is a *cross-extension of $H$ at $u$*. Assume now that $H$ is contained in a brace $G$. We say that $H$ is *$G$-flexible at $u$* if there exists a vertex $x \in A - \{u\}$ and a set $I \subseteq \{1,2,3\}$ of size two such that for every $i \in I$, $H+(u_i,x)$ is weakly contained in $G$.

5.1. *Let $H$ be a connected bipartite graph on at least four vertices, let $u$ be a vertex of $H$ of degree three, and let $G$ be a brace. If $G$ contains a cross-extension of $H$ at $u$, then $H$ is $G$-flexible at $u$.*

*Proof.* Let $(A, B)$ be a bipartition of $H$ with $u \in A$, and let $H', H'', u_1, u_2,$ $u_3, a_1, b_1, a_2, b_2$ be as in the definition of cross-extension. By 4.9 applied to $H''$ there exists a reduced $\{b_1, b_2\}$-augmentation $H'''$ of $H''$ that is weakly contained in $G$. From the symmetry we may assume that $H'''$ is determined by $b_1, e_1, e_2, \ldots, e_k$. From (R4) we deduce that $k = 1$. Thus $e_1 \in A - \{u\}$ by (R5). It follows that $H'''$ weakly contains $H + (u_1, e_1)$ and $H + (u_3, e_1)$, as desired. □

Let $H$ be a bipartite graph with bipartition $(A, B)$, let $u \in A$ have degree three, and let $u_1$, $u_2$ and $u_3$ be the neighbors of $u$. Let $H'$ be obtained from $H$ by replacing, for $i = 1, 2, 3$, the edge $uu_i$ by a path with vertex-set $\{u, b_i, a_i, u_i\}$ (in order), and let $H'' = H' + (b_1, a_2) + (b_2, a_3) + (b_3, a_1)$. We say that $H''$ is a *hexagonal extension of $H$ at $u$*.

5.2. *Let $H$ be a bipartite graph, let $u$ be a vertex of $H$ of degree three, and let $G$ be a brace. If $G$ contains a hexagonal extension of $H$ at $u$, then $H$ is $G$-flexible at $u$.*

*Proof.* Let $H'$, $H''$, $a_1$, $b_1$, $a_2$, $b_2$, $a_3$ and $b_3$ be as in the definition of hexagonal extension. Let $X = \{b_1, b_2, b_3\}$. By 4.9 some reduced $X$-augmentation $H'''$ of $H''$ is weakly contained in $G$. From the symmetry we may assume that $H'''$ is determined by $b_1, e_1, e_2, \ldots, e_k$. If $e_1 \in A$ (that is, $k = 1$), then $e_1$ satisfies the requirements for $H$ to be $G$-flexible at $u$. We may therefore assume that $k > 1$. By (R4) $e_1 \neq a_1 u_1$ and $e_1 \neq a_2 u_2$, and hence $e_1 = a_3 u_3$. It follows that $G$ contains a cross-extension of $H$ at $u$, and hence $H$ is $G$-flexible at $u$ by 5.1. □

5.3. *Let $H$ be a connected bipartite graph with maximum degree three, let $u$ be a vertex of $H$ of degree three, let $u_1, u_2, u_3$ be the neighbors of $u$ in $H$, and let $G$ be a brace. If $G$ weakly contains $H + (u_1, uu_2)$, then $H$ is $G$-flexible at $u$.*

*Proof.* Let $(A, B)$ be a bipartition of $H$ chosen so that $u \in A$. Let $H_1$ be obtained from $H$ by replacing the edge $uu_2$ by a path with vertices $u, b_2, a_2, u_2$ (in order). Then by 4.2 the graph $G$ contains $H_1' = H_1 + (a_2, y)$, where $y$ is $u_1$ or an edge of $H$ incident with $u_1$. If $y$ is an edge, let $b_1, a_1$ be the new vertices of $H_1'$; otherwise let $a_1 = b_1 = y$. In either case $G$ weakly contains a $b_2$-extension $H_2$ of $H_1'$ by 4.10. Suppose first that $H_2$ is of the first kind. Then $H_2 = H_1' + (b_2, x)$, where either $x \in A - \{u\}$, or $y$ is an edge and $x = a_1$. If $x \neq a_1$ then $x$ satisfies the requirements for $H$ to be $G$-flexible at $u$. This can be seen by considering the graph obtained from $H_2$ by deleting $ub_1$ or $uu_1$, and the graph obtained by deleting $b_1 a_2$. Thus we may assume that $y$ is an edge incident with $u_1$, and that $x = a_1$. If $y = uu_1$ and $x$ is as stated, then $H_2$ is a cross-extension of $H$ at $u$, and hence $H$ is $G$-flexible at $u$ by 5.1. Finally, if $y = vu_1$, where $v \neq u$ is a neighbor of $u_1$, then $H_2$ weakly contains $H + (u_2, v)$, and since $H = H + (u_1, v)$ we see that $v$ is as desired. Thus if $H_2$ is of the first kind the lemma holds.

Suppose now that $H_2$ is of the second kind. Let $H_2' = H_1' + (b_2, uu_3)$ and $H_2'' = H_1' + (b_2, a_2u_2)$, and in either case let $a_3, b_3$ be the new vertices. Then either $H_2 = H_2' + (b_3, x)$ or $H_2 = H_2'' + (b_3, x)$, where either $x \in A - \{u\}$, or $y$ is an edge and $x = a_1$. From the symmetry between $u$ and $a_2$ we may assume that $H_2 = H_2' + (b_3, x)$. If $y = uu_1$ and $x = a_1$ then $H_2$ is a hexagonal extension of $H$ at $u$, and thus $H$ is $G$-flexible at $u$ by 5.2. Thus we may assume that if $x = a_1$, then $y$ has ends $u_1$ and $x'$, where $x' \neq u$. If $x \neq a_1$ let $x' = x$. Now $H_2$ weakly contains $H + (x', u_1)$ (delete $a_2b_1$ and $ub_2$) and $H + (x', u_3)$ (delete $a_2b_1$ and $a_3b_2$), and hence $H$ is $G$-flexible at $u$.

Finally, suppose that $H_2$ is of the third kind, and let $H_2', H_2'', a_3, b_3$ be as in the previous paragraph. Again, from the symmetry we may assume that $H_2$ is $H_2' + (b_3, a_2u_2)$ or $H_2' + (b_3, a_2b_1)$. We claim that both of these graphs contain a cross-extension of $H$ at $u$. Indeed, in the first case delete $b_1a_2$, and in the second case delete $b_1a_4$, where $a_4, b_4$ are the new vertices of $H_2' + (b_3, a_2b_1)$. □

5.4. *Let $G$ be Rotunda with its vertex-set numbered as in Figure 3. Then $G + (5, 10)$ and $G + (5, 3)$ both contain $K_{3,3}$.*

*Proof.* To see the first containment delete the vertices 15 and 16, and the edges 2-10, 7-8 and 1-5. To see the second delete the vertices 1, 9, 13, 14, 15, 16. □

5.5. *Let $G$ be a brace not containing $K_{3,3}$, let $K$ be a model of Rotunda in $G$, and let $S$ be the set of four vertices of $K$ corresponding to the center of Rotunda. Then $G$ has no $(K, S)$-jump.*

*Proof.* Let $H$ be Rotunda with vertices numbered as in Figure 3. Let $H_1'$ be obtained from $H$ by replacing the edge 1-9 by a path with vertex set 1, 18, 17, 9 (in order), and let $H_1 = H_1' + (17, 2\text{-}6)$. Let us recall our earlier convention according to which the new vertices of $H_1$ are 19 and 20 in such a way that 20 has degree two.

Suppose for a contradiction that $G$ has a $(K, S)$-jump. Then $G$ weakly contains $H_1, H_2 = H + (5, 10)$, $H_3 = H + (5, 3)$, or $H_4 = H + (5, 1\text{-}9)$. But it does not weakly contain $H_2$ or $H_3$ by 5.4, and it does not weakly contain $H_4$ by 5.3, because it follows from 5.4 that $H$ is not $G$-flexible at 1. Thus $G$ weakly contains $H_1$, and hence by 4.2 it contains $H_1$.

By 4.10 $G$ weakly contains an 18-extension $H_5$ of $H_1$. Suppose first that $H_5$ is of the first kind. Then $H_5 = H_1 + (18, x)$, where $x \in \{3, 6, 8, 10, 12, 14, 16, 20\}$. If $x = 10$, then $H_5$ weakly contains $H_3$ (delete 9-10), if $x = 12$ then $H_5$ weakly contains $H_4$ (delete 9-12), if $x = 3$ then $H_5$ weakly contains $H_3$ (delete 17-19), and the same argument shows that in the remaining cases $H_5$ weakly contains $H_2$. Thus we may assume that $H_5$ is of the second or third kind. Since $H_1 + (18, 1\text{-}5)$ and $H_1 + (18, 1\text{-}13)$ weakly contain $H_4$ (delete 17-19), we may assume that $H_5 = H_1 + (18, 17\text{-}19) + (21, x)$, where $x \in \{3, 8, 14, 16, 1\text{-}5, 1\text{-}13\}$. It follows that $H_5$ weakly contains $H_2$, $H_3$ or $H_4$ (delete 17-22 and 19-21). Thus in each case $H_5$, and hence $G$, contains $K_{3,3}$, a contradiction. □

5.6. *Let $G$ be a brace, let $H$ be a subgraph of $G$ on at least four vertices, let $M$ be a perfect matching of $G \backslash V(H)$, and let $e \in M$ have ends $u$ and $v$. Then there exist two $M$-alternating paths $P_1$ and $P_2$ not containing $e$ with one end $u$, the other end in $H$ and with $V(P_1) \cap V(P_2) = \{u\}$.*

*Proof.* Let $(A, B)$ be a bipartition of $G$ with $u \in A$, and let $D$ be the digraph obtained from $G$ by directing every edge toward $B$ and contracting the edges of $M$. Since $G$ is 2-extendable it follows that $D$ has two directed paths from $u$ to $V(H)$, vertex-disjoint except for $u$. These paths give rise to the desired $M$-alternating paths in $G$ in the natural way. □

By a *walk* in a graph $G$ we mean a sequence $v_1, v_2, \ldots, v_n$ of vertices of $G$ such that $v_i$ and $v_{i+1}$ are adjacent for every $i = 1, 2, \ldots, n - 1$. We say that $v_1$ and $v_n$ are the *ends* of the walk. We say that a walk $v_1, v_2, \ldots, v_n$ is *M-alternating*, where $M$ is a matching in $G$, if the edge joining $v_i$ and $v_{i+1}$ belongs to $M$ for every $i \in I$, where $I$ is either the set of all odd or the set of all even integers in $\{1, 2, \ldots, n - 1\}$. We need the following easy but useful lemma.

5.7. *Let $M$ be a matching in a bipartite graph $G$, and let $u, v \in V(G)$. If there exists an $M$-alternating walk in $G$ with ends $u$ and $v$, then there exists an $M$-alternating path in $G$ with ends $u$ and $v$.*

Let $G$ be a bipartite graph with bipartition $(A, B)$. Let $v, v_1, v_2, v_3$ be distinct vertices of $G$ such that $v_1, v_2$ and $v_3$ belong to the same color class, and $v$ belongs to the opposite color class, and for $i = 1, 2, 3$ let $P_i$ be a path in $G$ with ends $v$ and $v_i$. If $P_1, P_2$ and $P_3$ pairwise intersect only in $v$ we say that $F = P_1 \cup P_2 \cup P_3$ is a *fork* in $G$. The vertices $v_1, v_2, v_3$ are called the *ends* of $F$.

Let $P$ be a path in $G$ with distinct ends $v_0$ and $v_{2k+1}$, and let $v_0, v_1, v_2, \ldots, v_{2k}, v_{2k+1}$ be some of the vertices of $P$ occurring on $P$ in the order listed and such that for $i = 0, 1, 2, \ldots, k$, $v_{2i}$ and $v_{2i+1}$ belong to opposite color classes, and for $i = 1, 2, \ldots, k$, $v_{2i-1}$ and $v_{2i}$ belong to the same color class (possibly $v_{2i-1} = v_{2i}$). For $i = 1, 2, \ldots, k$ let $F_i$ be a subgraph of $G$ with the following properties. If $v_{2i-1} = v_{2i}$, then $F_i$ is a fork with one end $v_{2i}$, and otherwise disjoint from $P$. If $v_{2i-1} \neq v_{2i}$, then $F_i = F_i' \cup F_i''$, where $F_i'$ and $F_i''$ are vertex-disjoint, $F_i'$ is a fork with two of its ends $v_{2i-1}$ and $v_{2i}$, and otherwise disjoint from $P$, and $F_i''$ is a path with ends $v \in V(P)$ and $w \notin V(P)$ and otherwise disjoint from $P$ in such a way that $v_{2i-1}, v, v_{2i}$ occur on $P$ in the order listed, and $v, v_{2i}$ belong to opposite color classes and $v, w$ belong to opposite color classes. The vertices of $F_i$ of degree one in $F_i$ are called the *ends* of $F_i$. Assume moreover that the graphs $F_i$ are pairwise disjoint except possibly for their ends. In those circumstances we call $K = P \cup F_1 \cup F_2 \cup \ldots \cup F_k$ a *hook*. We say that $K$ has *strength* $k$, that $v_0$ and $v_{2k+1}$ are the *ends* of $K$, and that the ends of the graphs $F_i$ that do not belong to $P$ are the *anchors* of $K$. Let $H$ be a subgraph of $G$. We say that a hook $K$ is an *H-hook* if $K$ has and only has its ends and its anchors in common with $H$. Thus an $H$-path with ends in opposite color classes is an $H$-hook of strength 0.

5.8. *Let $G$ be a brace, let $H$ be a subgraph of $G$ such that $G \backslash V(H)$ has a perfect matching, let $S \subseteq V(H)$, and assume that there exists an $H$-path in $G$ between different components of $H \backslash S$. Then there exists an $H$-hook $K$ in $G$ with ends in different components of $H \backslash S$ and with all anchors in $S$, such that $G \backslash V(H \cup K)$ has a perfect matching.*

*Proof.* Let $G, H, S$ be as stated, and let $M$ be a perfect matching of $G \backslash V(H)$. We proceed by induction on $|V(G)| - |V(H)|$. Let $P$ be an $H$-path in $G$ with ends in different components of $H \backslash S$; let one end of $P$ belong to $V(H_1)$, where $H_1$ is a component of $H \backslash S$. If $|V(P)| = 2$ then $P$ satisfies the conclusion of the theorem, and so we assume that $|V(P)| > 2$. Let $D$ be the set of all vertices of $G \backslash V(H)$ that can be reached from $V(H_1)$ by an $M$-alternating path. Then $D \neq \emptyset$ since $|V(P)| > 2$, and $|D \cap A| = |D \cap B|$, where $(A, B)$ is a bipartition of $G$. Let $H'$ be the subgraph of $G$ induced by $V(H) \cup D$; then $G \backslash V(H')$ has a perfect matching (namely a subset of $M$).

If there is no $H'$-path between different components of $H' \backslash S$ (for instance, this happens when $H' \backslash S$ is connected), then some vertex of $D$ is adjacent to a vertex in $V(H) - V(H_1) - S$, and otherwise we may apply the induction hypothesis to $G$, $H'$ and $S$. Thus in either case there exists an $H'$-hook $K'$ in $G$ with anchors in $S$, either with ends in different components of $H' \backslash S$, or one end in $D$ and the other end in $V(H) - V(H_1) - S$, and such that $G \backslash V(H' \cup K')$ has a perfect matching. If $K'$ has ends in different components of $H$, then $K'$ satisfies the conclusion of 5.8, and so we may assume that one end of $K'$, say $u$, belongs to $D$ and the other to $V(H) - V(H_1) - S$. Let $e \in M$ be the unique member of $M$ incident with $u$, and let $v$ be the other end of $e$. Let $Q$ be an $M$-alternating path with one end $u$ and the other end in $V(H_1)$. Then $V(Q) \cap V(K') = \{u\}$, because $V(Q) \subseteq D \cup V(H_1)$. If $v \in V(Q)$, then $Q \cup K'$ satisfies the conclusion of the lemma, and so we may assume that $v \notin V(Q)$. By 5.6 there exist two $M$-alternating paths $Q_1, Q_2$ not containing $e$ with one end $v$ and the other end in $V(H)$, disjoint except for $v$. Let $i \in \{1, 2\}$ and let $v_i \in V(H)$ be the other end of $Q_i$. Since $Q_i$ is $M$-alternating and $e \notin E(Q_i)$, it follows that $u \notin V(Q_i)$. Moreover, $Q_i$ is a subgraph of $H'$ by 5.7, because $v \in V(H')$. Thus $V(Q_i) \cap V(K') = \emptyset$. Let $L$ be the graph with $V(L) = \{u, v\}$ and $E(L) = \{e\}$. If $v_i \in V(H_1)$, then $K' \cup Q_i \cup L$ satisfies the conclusion of the lemma; if $v_i \in V(H) - V(H_1) - S$, then by 5.7 $Q_i \cup Q$ includes an $M$-alternating path between different components of $H \backslash S$, which satisfies the conclusion of the lemma. Thus we may assume that $v_1, v_2 \in S$. Then $L \cup Q_1 \cup Q_2$ is a fork with ends $u, v_1, v_2$ such that $G \backslash V(H \cup L \cup Q_1 \cup Q_2)$ has a perfect matching. Thus we have shown that given $G, H, S, M, H_1, K', u, v, e$ as above we may assume that there exist an $M$-alternating path $Q$, distinct vertices $v_1, v_2 \in S$ and a fork $F$ with ends $u, v_1, v_2$ such that $v \notin V(Q)$, one end of $Q$ is $u$ and the other end belongs to $H_1$, the three paths comprising $F$ are $M$-alternating, $V(F) \subseteq D \cup \{v_1, v_2\}$ and $e \in E(F)$. Let us choose $M, Q$ and $F$ satisfying these requirements and, subject to that, with $F \cup Q$ minimum.

Let $P$, $P_1$ and $P_2$ be the three paths comprising $F$ in such a way that $u$ is an end of $P$, $v_1$ is an end of $P_1$ and $v_2$ is an end of $P_2$, and let $w$ be the common end of $P$, $P_1$ and $P_2$. If $V(Q) \cap V(F) = \{u\}$, then $K' \cup Q \cup F$ is an $H$-hook satisfying the conclusion of the lemma, and so we may assume that $V(Q) \cap V(F) - \{u\} \neq \emptyset$. Let $x$ be the vertex of $V(Q) \cap V(F) - \{u\}$ chosen so that the subpath $Q'$ of $Q$ with one end $x$ and the other in $V(H)$ is as short as possible. Since $Q$, $P$, $P_1$ and $P_2$ are $M$-alternating and $e \notin E(Q)$ we deduce that $x$ and $u$ belong to the same color class. Thus if $x \in V(P)$, then $K' \cup Q' \cup F$ satisfies the conclusion of the theorem. By the symmetry between $P_1$ and $P_2$ we may therefore assume that $x \in V(P_1)$. Let $y \in V(Q)$ be chosen so that the subpath of $Q$ between $x$ and $y$ is a subgraph of $F$, and subject to that, the path has maximum length. Since $P_1$ is $M$-alternating and $w$ and $v_1$ belong to different color classes, it follows that $y$ belongs to the subpath of $P_1$ with ends $v_1$ and $x$. Thus $Q$ includes an $F$-path $Q''$ with one end $y$; let $z$ be the other end of $Q''$. Since $Q$ is $M$-alternating it follows that $z$ and $u$ belong to the same color class. If $z$ belongs to the subpath of $P_1$ between $y$ and $v_1$, then replacing the path of $F$ with ends $y$ and $z$ by $Q''$ results in a fork $F'$ that contradicts the minimality of $F \cup Q$. If $z$ belongs to the subpath of $P_1$ between $w$ and $y$, then let $F'$ be defined as above, let $C$ be the unique circuit of $F \cup Q''$, and let $M' = (M - E(C)) \cup (E(C) - M)$. Also, the union of $Q \backslash V(Q'')$ and the subpath of $P_1$ with ends $x$ and $z$ includes an $M'$-alternating walk with the same ends as $Q$, and hence it includes an $M'$-alternating path $Q_1$ with the same ends by 5.7. It follows that the triple $M', Q_1, F'$ contradicts the minimality of $F \cup Q$. Thus $z \notin V(P_1)$. If $z \in V(P_2) - \{w\}$, then $F \cup Q''$ includes a fork that contradicts the minimality of $F \cup Q$ (with $M$ replaced by a suitable matching). Thus $z \in V(P)$, and it follows that $K' \cup Q' \cup Q'' \cup F$ is a hook of strength one and that it satisfies the conclusion of the lemma. □

5.9. *Let $G$ be a brace, let $H$ be a model of Rotunda in $G$, and let $S$ be the set of vertices of $H$ corresponding to the center of Rotunda. If there exists an $H$-path in $G\backslash S$ with ends in different components of $H\backslash S$, then $G$ contains $K_{3,3}$.*

*Proof.* By 5.8 there exists an $H$-hook $K$ in $G$ with ends in different components of $H\backslash S$ and with anchors in $S$, and such that $G\backslash V(H\cup K)$ has a perfect matching. From 5.5 $K$ has strength at least one. Let us assume the notation introduced prior to 5.8. If $K$ has strength at least two, then using $F_1, F_2$ and the subpath of $P$ between them it is easy to establish that $G$ contains $K_{3,3}$. (To see this notice that the anchors of $K$ that belong to $F_1$ and $F_2$ are distinct, and hence are equal to $1, 2, 3, 4$.) Thus we may assume that $K$ has strength one. Let $R$ be Rotunda with vertices numbered as in Figure 3. From the symmetry we may assume that the anchors of $K$ are the vertices of $H$ that correspond to the vertices 1 and 3 of $R$. We deduce (using symmetry of $R$) that $G$ weakly contains $R + (3,5) + (1, 3\text{-}5) + (18, 1\text{-}17) + (20, x)$, or $R + (3,5) + (1, 3\text{-}5) + (18, x)$, or $G$ contains $R_1 = R + (3, 2\text{-}6) + (1, 3\text{-}17) + (20, 1\text{-}19) + (22, f)$, or $R_2 = R + (3, 2\text{-}6) + (1, 3\text{-}17) + (20, f)$, where $x$ is a suitable vertex or edge of $R$, and $f$ is 2-10 or 4-12. (Notice that if $f$ is the edge 1-9, then by bicontracting we may assume that $f = 9$, which is a case already covered, and similarly if $f$ is the edge 3-11.) The first two graphs contain $K_{3,3}$ by 5.4, and so we may assume that $G$ contains $R_1$ or $R_2$.

We assume first that $G$ contains $R_1$. By 4.10 $G$ weakly contains an 18-extension $R'_1$ of $R_1$. By 5.5 we may assume that $G$ has no $(H, S)$-jump, and so by the symmetry of $R$ we may assume that $R'_1 = R_1 + (7, 18)$, $R'_1 = R_1 + (4, 18)$, $R'_1 = R_1 + (18, 19)$, or $R'_1 = R_1 + (18, 17\text{-}20) + (4, 26)$. The first two and the fourth graph contain $K_{3,3}$ by 5.4. In more detail, the first weakly contains $R + (3, 5)$ by deleting 21, 22 and 6-7, the second weakly contains $R + (4, 6)$ by deleting 19, 20, 21 and 22, and the fourth weakly contains $R + (4, 6)$ by deleting 19, 20, 21, 22 and 17-26. If $G$ weakly contains the third graph, then $G$ contains $K_{3,3}$ by 5.5 (delete the vertices 20 and 21 and the edge 3-19). Thus if $G$ contains $R_1$ then $G$ contains $K_{3,3}$.

We may therefore assume that $G$ contains $R_2$. By 4.10 $G$ weakly contains an 18-extension $R'_2$ of $R_2$. By 5.5 we may assume that $G$ has no $(H, S)$-jump, and so by the symmetry of $R$ we may assume that $R'_2 = R_2 + (7, 18)$, $R'_2 = R_2 + (18, 4)$, $R'_2 = R_2 + (18, 19)$, or $R_2 = R_1 + (18, 17\text{-}20) + (4, 24)$. The first graph weakly contains $R + (3, 5)$ (delete 6-7, 20-21 and 1-19), and the second and fourth graphs weakly contain $R + (6, 4)$. In those cases $G$ contains $K_{3,3}$ by 5.4 and the symmetry of Rotunda. Finally, if $G$ weakly contains the third graph, then $G$ contains $K_{3,3}$ by 5.5 (delete 1-19, 3-19 and 17-20). □

## 6. Proof of the main result

We start with six lemmas. The first is a theorem of Kasteleyn [6].

6.1. *Every planar graph admits a Pfaffian orientation.*

6.2. *The graph $K_{3,3}$ does not admit a Pfaffian orientation.*

*Proof.* Let $(A, B)$ be a bipartition of $K_{3,3}$, and let $\mathcal{C}$ be the set of all circuits of $K_{3,3}$ of length four. Then $|\mathcal{C}| = 9$. Suppose for a contradiction that $D$ is a Pfaffian orientation of $K_{3,3}$, let $\varepsilon$ be the number of edges directed in $D$ from $A$ to $B$, and for $C \in \mathcal{C}$ let $\varepsilon_C$ be the number of edges of $C$ directed in $D$ from $A$ to $B$. Then $\varepsilon_C$ is odd, because $C$ is central. Since every edge of $K_{3,3}$ belongs to four members of $\mathcal{C}$, we have $4\varepsilon = \sum_{C \in \mathcal{C}} \varepsilon_C$, a contradiction, because the right-hand side is odd. □

6.3. *The Heawood graph admits a Pfaffian orientation.*

*Proof.* Let $H$ be the Heawood graph and let $(A, B)$ be a bipartition of $H$. Orienting every edge of $H$ from $A$ to $B$ gives a Pfaffian orientation of $H$, because no circuit of $H$ of length eight or twelve is central. □

6.4. *Let a bipartite graph $G$ contain a bipartite graph $H$. If $G$ admits a Pfaffian orientation then so does $H$.*

*Proof.* Let $D$ be a Pfaffian orientation of $G$, and let $\phi$ be an embedding of $H$ into $G$. Let $e$ be an edge of $H$ with ends $u$ and $v$. Then $\phi(e)$ is a path with an odd number of edges; if an odd number of those edges are directed from $\phi(u)$ to $\phi(v)$ we direct $e$ from $u$ to $v$, and otherwise we direct $e$ from $v$ to $u$. This defines a Pfaffian orientation of $H$, because the image under $\phi$ of a central circuit in $H$ is a central circuit in $G$. □

If $G$, $G_1$, $G_2$ and $C$ are as in the definition of 4-sum, we say that $G$ *is a 4-sum of $G_1$ and $G_2$ along $C$*.

6.5. *Let $G$ be a bipartite graph such that $G$ is a 4-sum of $G_1$ and $G_2$ along $C$.*

(i) *If $M$ is a perfect matching of $G$, then $E(G_1) \cap M$ is a subset of a perfect matching of $G_1$.*

(ii) *If $G_1$ and $G_2$ admit Pfaffian orientations, then so does $G$.*

(iii) *If $G$ is a brace, then so are $G_1$ and $G_2$.*

*Proof.* To prove (i) let $M$ be a perfect matching of $G$. Let us say that a vertex $v \in V(C)$ is *exposed* if the unique edge of $M$ incident with $v$ does not belong to $G_1$. Since $C$ is central we deduce that $C$ has a matching $M'$ such that the set of exposed vertices of $G_1$ is precisely the set of vertices incident with an edge in $M'$. It follows that $(E(G_1) \cap M) \cup M'$ is a perfect matching in $G_1$. This proves (i).

If $D$ is a Pfaffian orientation of a graph $G$, $v$ is a vertex of $G$, and $D'$ is obtained from $D$ by reversing the directions of all edges incident with $v$, then $D'$ is a Pfaffian orientation of $G$. Thus to prove (ii) we may choose Pfaffian orientations $D_1$ and $D_2$ of $G_1$ and $G_2$, respectively, such that the two orientations agree on $C$. Let $D$ be an orientation of $G$ defined by giving an edge of $G_i$ its orientation in $D_i$ ($i = 1, 2$). We claim that $D$ is a Pfaffian orientation of $G$. To prove this let $C_0$ be a central circuit in $G$; we must show that $C_0$ is oddly oriented in $D$. We proceed by induction on $|V(C_0)|$. If $V(C_0) \subseteq V(G_1)$, let $M_1$ be a perfect matching of $C_0$ and let $M_2$ be a perfect matching of $G \backslash V(C_0)$. By (i) applied to $M_1 \cup M_2$ we see that $C_0$ is a central circuit in $G_1$, and so is oddly oriented in $D_1$, and hence it is oddly oriented in $D$. A similar argument works when $V(C_0) \subseteq V(G_2)$, and so we may assume that $V(C_0) - V(G_1) \neq \emptyset \neq V(C_0) - V(G_2)$. Assume first that it is not the case that $V(C) \cap V(C_0)$ consists of two diagonally opposite vertices of $C$. Then there exist circuits $C_1$ and $C_2$ of $G$ such that $C_0 = C_1 \cup C_2 \backslash e$ for some edge $e$ of $C$. For $i = 1, 2$ the graph $G \backslash V(C_i)$ has a perfect matching (namely the union of a perfect matching of $G \backslash V(C_0)$ and a perfect matching of $V(C_{3-i}) \backslash V(C_i)$), and hence both $C_1$ and $C_2$ are oddly oriented in $D$ by the induction hypothesis. It follows that $C_0$ is oddly oriented, as desired. We may therefore assume that $V(C) \cap V(C_0) = \{u_1, u_3\}$, where $u_1, u_2, u_3, u_4$ are the vertices of $C$ in order. As before, let $M_2$ be a perfect matching of $G \backslash V(C_0)$. We may assume without loss of generality that the edge of $M_2$ incident with $u_2$ belongs to $G_1$, and hence the edge of $M_2$ incident with $u_4$ belongs to $G_2$. Let $P_1, P_2$ be the two subpaths of $C_0$ with ends $u_1$ and $u_3$ and union $C_0$, and let $Q_1, Q_2$ be defined similarly with $C$ replacing $C_0$. We may assume that $u_4 \in V(Q_1)$, $u_2 \in V(Q_2)$, and that $P_i$ is a subgraph of $G_i$ for $i = 1, 2$. Then the circuits $P_1 \cup Q_1$ and $P_2 \cup Q_2$ are central in $G_1$ and $G_2$, respectively, and hence they are both oddly oriented in $D$ by the induction hypothesis. It follows that $C_0$ is oddly oriented, as desired. This proves (ii).

To prove (iii) let us assume that $G$ is a brace. By 3.1 every bipartite graph obtained from $G$ by adding edges is a brace, and so we may assume that $C$ is a subgraph of $G$. Let $e, f$ be edges of $G_1$ with no common end; then $e, f \in E(G)$, and so there exists a perfect matching $M$ of $G$ containing $e$ and $f$. By (i) $E(G_1) \cap M$ is a subset of a perfect matching $M'$ of $G_1$; then $e, f \in M'$, as desired. □

6.6. *Let $H$ be a connected graph contained in a connected brace $G$. Let every vertex of $H$ have degree three and let $H$ have no circuit of length four. Then either $H$ is isomorphic to $G$, or there exist nonadjacent vertices $u$ and $v$ of $H$ belonging to different color classes of $H$ such that $H + (u, v)$ is weakly contained in $G$.*

*Proof.* Let us assume that $H$ is not isomorphic to $G$, and let $K$ be a model of $H$ in $G$ with $|V(K)|$ maximum. Assume first that $K$ has no vertex of degree two; then $K \neq G$. Let $M$ be a perfect matching of $G \backslash V(K)$. There exists an edge $e \in E(G) - E(K)$ incident with a vertex of $K$. Let $M'$ be a perfect matching of $G$ containing $e$. Then $M \triangle M'$ includes the edge-set of a $K$-path $P$ such that $G \backslash V(K \cup P)$ has a perfect matching. Let $u'$ and $v'$ be the ends of $P$; since $K$ has no vertex of degree two the vertices $u'$ and $v'$ correspond to vertices $u$ and $v$ of $H$, respectively. If $u$ and $v$ are not adjacent in $H$, then $G$ contains $H + (u, v)$, and the lemma holds. Thus we may assume that $u$ and $v$ are adjacent in $H$, and hence $u'$ and $v'$ are adjacent in $K$. Let $K' = (K \backslash u'v') \cup P$; then $K'$ is a model of $H$ in $G$, and hence contradicts the choice of $K$.

We may therefore assume that $K$ has a vertex of degree two. Thus there exists an edge $e$ of $H$ such that the corresponding path of $K$ has at least two internal vertices. Let $(A, B)$ be a bipartition of $H$, let the ends of $e$ be $u \in A$ and $v \in B$, and let $H_1$ be obtained from $H$ by replacing the edge $e$ by a path with vertices $u, v', u', v$ (in order). By 4.10 $G$ contains either $H_1 + (u', w)$ or $H_1 + (u', e)$, where $w \in B - \{v\}$ and $e$ is an edge of $H_1$ incident with $v$ but not with $u'$. We deduce that either the lemma holds, or $G$ weakly contains $H + (x, uv)$ or $H + (u, vy)$, where $x \neq v$ is a neighbor of $u$ in $H$ and $y \neq u$ is a neighbor of $v$ in $H$. In the last two cases $H$ is $G$-flexible at $u$ or $v$ by 5.3. From the symmetry we may assume that $H$ is $G$-flexible at $u$; that is, there exists a vertex $x \in V(H) - \{u\}$ and two distinct neighbors $u_1$ and $u_2$ of $u$ such that $x$ and $u$ belong to the same color class and $G$ weakly contains both $H + (u_1, x)$ and $H + (u_2, x)$. Since $H$ has no circuit of length four it follows that $x$ is not adjacent to both $u_1$ and $u_2$. Say it is not adjacent to $u_1$; then $u_1$ and $x$ satisfy the conclusion of the lemma. □

6.7. *Let $G$ be a brace containing the Heawood graph and not containing $K_{3,3}$. Then $G$ is isomorphic to the Heawood graph.*

*Proof.* Let $H$ be the Heawood graph. If $G$ is not isomorphic to $H$, then by 6.6 it weakly contains $H + (u, v)$ for some pair $u, v$ of nonadjacent vertices of $H$ belonging to different color classes. But each such graph contains $K_{3,3}$, as is easily verified (in fact, they are all isomorphic, and so it suffices to check one graph), and hence so does $G$, a contradiction. □

Let $G_0$ be a graph, let $C$ be a central circuit of $G_0$ of length four, and let $G_1, G_2, G_3$ be three subgraphs of $G_0$ such that $G_1 \cup G_2 \cup G_3 = G_0$, and for distinct integers $i, j \in \{1, 2, 3\}$, $G_i \cap G_j = C$ and $V(G_i) - V(C) \neq \emptyset$. Let $G$ be obtained from $G_0$ by deleting some (possibly none) of the edges of $C$. In these circumstances we say that $G$ is a *trisum* of $G_1$, $G_2$ and $G_3$. We are now ready to prove our main theorem, which we restate in a slightly stronger form. The equivalence of (i) and (ii) is due to Little [9].

6.8. *For a brace $G$ the following conditions are equivalent*:

(i) *$G$ does not contain $K_{3,3}$.*
(ii) *$G$ has a Pfaffian orientation.*
(iii) *Either $G$ is isomorphic to the Heawood graph, or $G$ can be obtained from planar braces by repeated application of the 4-sum operation.*
(iv) *Either $G$ is isomorphic to the Heawood graph, or $G$ can be obtained from planar braces by repeated application of the trisum operation.*

*Proof.* From 6.4 and 6.2 we deduce that (ii) implies (i), from 6.1, 6.3 and 6.5(ii) we deduce that (iii) implies (ii), and clearly (iv) implies (iii). To prove that (i) implies (iv) let $G$ be a brace not containing $K_{3,3}$. We may assume that $G$ is not isomorphic to the Heawood graph.

We prove by induction on $|V(G)|$ that $G$ can be obtained by repeated application of the trisum operation as stated in the theorem. If $G$ is planar then the claim holds, and so we may assume that $G$ is not planar. By 1.5 $G$ contains the Heawood graph or Rotunda. By 6.7 it does not contain the Heawood graph, and hence it contains Rotunda. Let $H$ be a model of Rotunda in $G$, let $S$ be the set of four vertices of $H$ that correspond to the center of Rotunda, and let $G_0$ be obtained from $G$ by adding an edge with ends $u$ and $v$, for every pair of vertices $u, v \in S$ that are not adjacent in $G$ and belong to different color classes of $G$. Thus $S$ induces a circuit $C$ in $G_0$. Moreover, $C$ is central in $G_0$, because $G$ is a brace. By 5.9 there exist graphs $G_1$, $G_2$ and $G_3$ such that $G_1 \cup G_2 \cup G_3 = G_0$, for distinct integers $i, j \in \{1, 2, 3\}$ $G_i \cap G_j = C$ and for $i = 1, 2, 3$ exactly one component of $H \backslash S$ is a subgraph of $G_i$. We claim that each $G_i$ is contained in $G$. From the symmetry it suffices to argue for $i = 1$. We have $E(G_1) - E(C) \subseteq E(G)$, and so it remains to account for the edges of $C$. One perfect matching of $C$ may be represented in $G_2$ (using two disjoint paths of $H \backslash S$) and the complementary matching may be represented in $G_3$ similarly. This proves that each $G_i$ is contained in $G$. Thus no $G_i$ contains $K_{3,3}$, and by 6.5(iii) each $G_i$ is a brace. Since the Heawood graph has no circuit of length four, it follows from the induction hypothesis that $G_1, G_2, G_3$ can be obtained from planar braces by repeated application of the trisum operation. Since $G$ is a trisum of $G_1$, $G_2$ and $G_3$, the theorem follows. □

## 7. Applications

Let $G$ be a bipartite graph with a perfect matching, let $(A, B)$ be a bipartition of $G$, and let $X$ be a nonempty proper subset of $A$ such that $|N(X)| = |X|$. Let $G_2 = G\backslash(X \cup N(X))$ and $G_1 = G\backslash V(G_2)$. We say that $G$ is a 0-*sum* of $G_1$ and $G_2$. We define 2-sum as follows. Let $G$ be a connected 1-extendable bipartite graph with bipartition $(A, B)$ that is not 2-extendable. Then $G$ has an edge $e$ with ends $u_1 \in A$ and $u_2 \in B$ such that $G\backslash u_1 \backslash u_2$ is not 1-extendable. By 3.1 applied to every component of $G\backslash u_1\backslash u_2$ there exists a nonempty proper subset $X$ of $A - \{u_1\}$ such that $|Y_1| = |X|$, where $Y_1 = N_G(X) - \{u_2\}$. Let $Y_2 = A - X - \{u_1\}$, let $G_2' = G\backslash(X \cup Y_1)$ and let $G_1' = G\backslash(V(G_2') - \{u_1, u_2\})$. For $i = 1, 2$ let $Y_i'$ be the set of all vertices of $Y_i$ that are not adjacent to $u_i$ but are adjacent to a vertex of $G\backslash V(G_i')$, and let $G_i$ be obtained from $G_i'$ by joining each vertex of $Y_i'$ by an edge to $u_i$. We say that $G$ is a 2-*sum* of $G_1$ and $G_2$. The following lemma follows from [10].

7.1. *Let $G_1$ and $G_2$ be bipartite graphs, let $i \in \{0, 2\}$, and let $G$ be an $i$-sum of $G_1$ and $G_2$. Then $G$ has a Pfaffian orientation if and only if both $G_1$ and $G_2$ have Pfaffian orientations.*

We deduce the following refinement of our main theorem.

7.2. *A bipartite graph has a Pfaffian orientation if and only if it either has no perfect matching, or it can be obtained by repeatedly applying the 0-, 2- and 4-sum operations, starting from connected planar bipartite graphs with perfect matchings and the Heawood graph.*

As a corollary of our main theorem we get the following extremal result.

7.3. *Every brace with $n \geq 3$ vertices and more than $2n - 4$ edges contains $K_{3,3}$, and hence does not have a Pfaffian orientation.*

*Proof.* Every planar bipartite graph on $n \geq 3$ vertices has at most $2n - 4$ edges. The result follows from 1.3 by induction. □

Let us turn to directed graphs now. A *directed graph* $D$ (or *digraph* for short) consists of a finite set $V(D)$ of vertices, a finite set $E(D)$ of edges, and an incidence relation that assigns to each edge of $D$ an ordered pair of distinct vertices of $D$ in such a way that different edges are assigned different ordered pairs. If $(u, v)$ is the ordered pair assigned to the edge $e$, we say that $u$ is the *tail* of $e$, and that $v$ is the *head* of $e$, and we denote the edge $e$ by $uv$. *Circuits* in digraphs are *directed*, and have no "repeated" vertices. A digraph $D$ is *even* if for every weight function $w : E(D) \to \{0, 1\}$ there exists a circuit in $D$ of even total weight. It was shown in [19] and is not difficult to see that testing evenness is polynomial-time equivalent to testing whether a digraph has an even directed circuit. Let $G$ be a bipartite graph with bipartition $(A, B)$, and let $M$ be a perfect matching in $G$. Let $D = D(G, M)$ be obtained from $G$ by directing every edge from $A$ to $B$, and contracting every edge of $M$. Little [9] has shown the following.

7.4. *Let $G$ be a bipartite graph, and let $M$ be a perfect matching in $G$. Then $G$ has a Pfaffian orientation if and only if $D(G, M)$ is not even.*

Since every digraph is isomorphic to $D(G, M)$ for some $G$ and $M$, 7.2 gives a characterization of even directed graphs. Let us state the characterization explicitly, but first let us point out a relation between extendability and strong connectivity. A digraph $D$ is *strongly connected* if for every two vertices $u$ and $v$ it has a directed path from $u$ to $v$. It is *strongly $k$-connected*, where $k \geq 1$ is an integer, if for every set $X \subseteq V(D)$ of size less than $k$, the digraph $D \backslash X$ is strongly connected. The following is straightforward.

7.5. *Let $G$ be a connected bipartite graph, let $M$ be a perfect matching in $G$, and let $k \geq 1$ be an integer. Then $G$ is $k$-extendable if and only if $D(G, M)$ is strongly $k$-connected.*

Let $D$ be a digraph, and let $(X, Y)$ be a partition of $V(G)$ into two nonempty sets in such a way that no edge of $G$ has tail in $X$ and head in $Y$. Let $D_1 = D \backslash Y$ and $D_2 = D \backslash X$. We say that $D$ is a *0-sum* of $D_1$ and $D_2$. Now let $v \in V(D)$, and let $(X, Y)$ be a partition of $V(D) - \{v\}$ into two nonempty sets such that no edge of $D$ has tail in $X$ and head in $Y$. Let $D_1$ be obtained from $D$ by deleting all edges with both ends in $Y \cup \{v\}$ and identifying all vertices of $Y \cup \{v\}$, and let $D_2$ be obtained by deleting all edges with both ends in $X \cup \{v\}$ and identifying all vertices of $X \cup \{v\}$. We say that $D$ is a *1-sum* of $D_1$ and $D_2$. Let $D_0$ be a directed graph, let $u, v \in V(D_0)$, and let $uv, vu \in E(D_0)$. Let $D_1$ and $D_2$ be such that $D_1 \cup D_2 = D_0$, $V(D_1) \cap V(D_2) = \{u, v\}$, $V(D_1) - V(D_2) \neq \emptyset \neq V(D_2) - V(D_1)$ and $E(D_1) \cap E(D_2) = \{uv, vu\}$. Let $D$ be obtained from $D_0$ by deleting some (possibly neither) of the edges $uv, vu$. We say that $D$ is a *2-sum* of $D_1$ and $D_2$. Now let $D_0$ be a directed graph, let $u, v, w \in V(D_0)$, let $uv, wv, wu \in E(D_0)$, and assume that $D_0$ has a directed circuit containing the edge $wv$, but not the vertex $u$. Let $D_1$ and $D_2'$ be such that $D_1 \cup D_2' = D_0$, $V(D_1) \cap V(D_2') = \{u, v, w\}$, $V(D_1) - V(D_2') \neq \emptyset \neq V(D_2') - V(D_1)$ and $E(D_1) \cap E(D_2') = \{uv, wv, wu\}$, let $D_2'$ have no edge with tail $v$, and no edge with head $w$. Let $D$ be obtained from $D_0$ by deleting some (possibly none) of the edges $uv, wv, wu$, and let $D_2$ be obtained from $D_2'$ by contracting the edge $wv$. We say that $D$ is a *3-sum* of $D_1$ and $D_2$. Finally let $D_0$ be a directed graph, let $x, y, u, v \in V(D_0)$, let $xy, xv, uy, uv \in E(D_0)$, and assume that $D_0$ has a directed circuit containing precisely two of the edges $xy, xv, uy, uv$.

Let $D_1$ and $D_2'$ be such that $D_1 \cup D_2' = D_0$, $V(D_1) \cap V(D_2') = \{x, y, u, v\}$, $V(D_1) - V(D_2') \neq \emptyset \neq V(D_2') - V(D_1)$ and $E(D_1) \cap E(D_2') = \{xy, xv, uy, uv\}$, let $D_2'$ have no edge with tail $y$ or $v$, and no edge with head $x$ or $u$. Let $D$ be obtained from $D_0$ by deleting some (possibly none) of the edges $xy, xv, uy, uv$, and let $D_2$ be obtained from $D_2'$ by contracting the edges $xy$ and $uv$. We say that $D$ is a 4-*sum* of $D_1$ and $D_2$. We say that a digraph is *strongly planar* if it has a planar drawing such that for every vertex $v \in V(D)$, the edges of $D$ with head $v$ form an interval in the cyclic ordering of edges incident with $v$ determined by the planar drawing. Let $F_7$ be the directed graph $D(H, M)$, where $H$ is the Heawood graph, and $M$ is a perfect matching of $H$. This defines $F_7$ uniquely up to isomorphism, irrespective of the choice of the bipartition of $H$ or the choice of $M$. Theorems 6.8 and 7.4 imply the following.

7.6.   *A digraph $D$ is not even if and only if it can be obtained from strongly planar digraphs and $F_7$ by means of 0-, 1-, 2-, 3- and 4-sums.*

The proof is fairly straightforward, and we omit it. Let us point out, however, that it requires the equivalence of 6.8(ii) and 6.8(iv), as opposed to merely 1.3.

From 7.3, 7.4 and 7.5 we deduce the following extremal result.

7.7.   *Let $D$ be a strongly 2-connected directed graph on $n \geq 2$ vertices. If $D$ has more than $3n - 4$ edges, then $D$ is even.*

Corollary 7.7 does not hold for strongly connected digraphs. However, Thomassen [22] has shown that every strongly connected directed graph with minimum in- and out-degree at least three is even. This is equivalent to the following.

7.8.   *Let $G$ be a 1-extendable bipartite graph such that every vertex has degree at least four. Then $G$ contains $K_{3,3}$, and hence does not have Pfaffian orientation.*

*Proof.* Let $G$ be as stated. We may assume that $G$ is connected. Since $G$ has at least $2|V(G)|$ edges, we see that if $G$ is a brace the corollary follows from 7.3. We may therefore assume that $G$ is a 2-sum of $G_1$ and $G_2$, and we may choose $G_1$ and $G_2$ so that $|V(G_1)|$ is minimum. Clearly $|V(G_1)| \geq 4$. We claim that $G_1$ is a brace. Indeed, suppose for a contradiction that $G_1$ is not a brace. Let $(A, B)$ be a bipartition of $G$, and let $u_1 \in A$, $u_2 \in B$, $X \subseteq A - \{u_1\}$, $Y_1$ and $Y_2$ be as in the definition of 2-sum. Since $G_1$ is not a brace, it has a nonempty set $X' \subseteq A \cap V(G_1)$ such that $|N_{G_1}(X')| \leq |X'| + 1$ and $N_{G_1}(X') \neq B \cap V(G_1)$. Let $Y' = B - N_{G_1}(X')$. Then $|N_{G_1}(Y')| \leq |Y'| + 1$ and $N_{G_1}(Y') \neq A \cap V(G_1)$, and either $u_1 \notin X'$ or $u_2 \notin Y'$. Using $X'$ (if $u_1 \notin X'$) or $Y'$ (if $u_2 \notin Y'$) the graph $G$ can be expressed as a 2-sum of $G'_1$ and $G'_2$, where $|V(G'_1)| < |V(G_1)|$, contrary to the choice of $G_1$. This proves that $G_1$ is a brace. Since all vertices in $X$ have degree in $G_1$ at least four, and $u_1$ has degree in $G_1$ at least two, we see that $|E(G_1)| \geq 2|V(G_1)| - 2$, and so $G_1$ contains $K_{3,3}$ by 7.3. Since $G$ weakly contains $G_1$ we deduce from 4.2 that $G$ contains $K_{3,3}$, as desired. □

## 8. Toward an algorithm

Let $G$ be a brace, and let $G$ be a trisum of $G_1$, $G_2$ and $G_3$. By 6.5(ii), if $G_1$, $G_2$ and $G_3$ have Pfaffian orientations, then so does $G$. We now prove the converse. By Cube we mean the graph of (the 1-skeleton of) the 3-dimensional cube. Let $G$ be a graph, and let $C$ be a circuit in $G$ of length four. We say that $G$ is $C$-fat if $C$ is a circuit of some model of Cube in $G$. We say that $G$ has a $C$-cross if $C$ is a circuit of some model of $K_{3,3}$ in $G$. Our first objective is the following.

8.1. *Let $G$ be a connected brace, and let $C$ be a circuit in $G$ of length four. If $C \neq G$, then either $G$ is $C$-fat, or $G$ has a $C$-cross.*

*Proof.* Let $(A, B)$ be a bipartition of $G$, and let $u_1, u_2, u_3, u_4$ be the vertices of $C$ in order. We may assume that $u_1, u_3 \in A$ and $u_2, u_4 \in B$. Since $G$ is a brace, $G \backslash V(C)$ has a perfect matching $M$, and by 3.2 $u_1$ is incident with at least one edge $e \notin E(C)$. Let $M'$ be a perfect matching of $G$ containing $e$, and let $P_1$ be the component of the subgraph of $G$ with vertex-set $V(G)$ and edge-set $M \triangle M'$ that contains $e$. It follows that $P_1$ is a path with one end $u_1$ and the other end $u_2$ or $u_4$. From the symmetry we may assume that the other end is $u_2$. Moreover, $G \backslash V(C \cup P_1)$ has a perfect matching, namely $M - E(P_1)$. Let $X = V(P_1) \cap A - \{u_1\}$. Then $|N_{C \cup P_1}(X)| \leq |X| + 1$ and $N_{C \cup P_1}(X) \neq B \cap V(C \cup P_1)$. By 4.1 there exists a $(C \cup P_1)$-path $P_2$ in $G$ with one end $v_2 \in X$, the other end $u_4$ and such that $G \backslash V(C \cup P_1 \cup P_2)$ has a perfect matching. Let $P'_1$ be the subpath of $P_1$ with ends $v_2$ and $u_2$. We may assume that $P_1$ and $P_2$ are chosen so that $|V(P'_1)|$ is minimum. Let $Y = V(P_1) \cap B - V(P'_1)$; then $|N_{C \cup P_1 \cup P_2}(Y)| \leq |Y| + 1$ and $N_{C \cup P_1 \cup P_2}(Y) \neq A \cap V(C \cup P_1 \cup P_2)$. By 4.1 there exists a $(C \cup P_1 \cup P_2)$-path $P_3$ with one end $v_1 \in Y$ and the other end $v_4 \in A - N_{C \cup P_1 \cup P_2}(Y)$ such that $G \backslash V(C \cup P_1 \cup P_2 \cup P_3)$ has a perfect matching. If $v_4 = u_3$, then $C \cup P_1 \cup P_2 \cup P_3$ is isomorphic to $K_{3,3}$, and hence $G$ has a $C$-cross. If $v_4 \in V(P'_1)$, then the graph obtained from $P_1 \cup P_2 \cup P_3$ by deleting the interior of the subpath of $P_1$ with ends $v_1$ and $v_2$ contradicts the choice of $P_1$ and $P_2$. We may therefore assume that $v_4 \in A \cap V(P_2) - \{v_2\}$.

The graph $C \cup P_1 \cup P_2 \cup P_3$ proves that there exists a graph $H = C \cup C' \cup Q_1 \cup Q_2 \cup Q_4$ such that

(1) $G\backslash V(H)$ has a perfect matching,

(2) $C'$ is a circuit of $G$ disjoint from $C$, and

(3) for $i = 1, 2, 4$, $Q_i$ is a path in $G$ with an odd number of edges with one end $u_i$, the other end say $v_i \in V(C')$, and otherwise disjoint from $C \cup C'$.

We may assume that $H$ is chosen so that, subject to (1), (2) and (3),

(4) $Q_2 \cup Q_4$ is minimal.

Let $Z = V(C' \cup Q_1) \cap B$. Then $|N_H(Z)| \leq |Z| + 1$, and hence by 4.1 there exists an $H$-path $P$ with one end $u \in Z$ and the other $v \in A - N_H(Z)$ such that $G\backslash V(H \cup P)$ has a perfect matching. If $v \in A \cap V(Q_2 \cup Q_4)$, then by deleting a suitable subpath of $C'$ from $H \cup P$ we obtain a graph contradicting (4). Thus $v \notin V(Q_2 \cup Q_4)$, and hence $v = u_3$ (notice that $u_1 \in N_H(Z)$). Let $P_0$ be the subpath of $C'$ with ends $v_2, v_4$ that contains $v_1$. If $u \notin V(P_0)$, then $H \cup P$ is isomorphic to an even subdivision of Cube, and hence $G$ is $C$-fat. Thus we may assume that $u \in V(P_0 \cup Q_1)$. Then $H \cup P$ contains $K_{3,3}$ (delete the subpath of $C'$ with ends $v_1, v_4$ or $v_1, v_2$, depending on which does not include $u$), and hence $G$ has a $C$-cross. □

8.2. *Let $G$ be a connected brace that has a Pfaffian orientation, and let $G$ be a trisum of $G_1$, $G_2$ and $G_3$ along $C$. Then for $i = 1, 2, 3$, the graph $G_i$ is $C$-fat and does not have a $C$-cross.*

*Proof.* By 6.5(iii) each $G_i$ is a brace, and by 8.1 it is $C$-fat or has a $C$-cross. If some $G_i$ has a $C$-cross, then $G$ contains one of the graphs depicted in Figure 5. However, each of those graphs contains $K_{3,3}$ (the edges to be deleted are drawn thicker), contrary to 6.2 and 6.4. Thus no $G_i$ has a $C$-cross, and hence is $C$-fat by 8.1. □

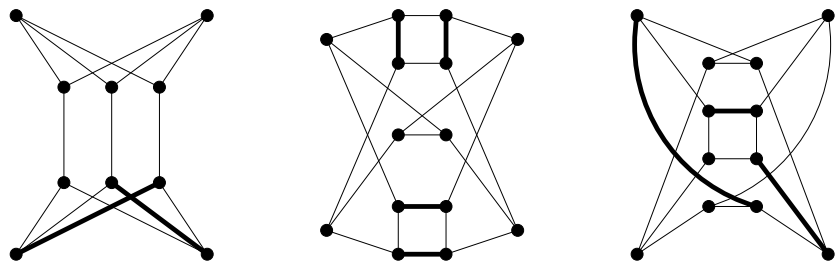

Figure 5

8.3. *Let $G$ be a connected brace that has a Pfaffian orientation, and let $G$ be a trisum of $G_1$, $G_2$ and $G_3$ along $C$. Then $G_1$, $G_2$ and $G_3$ have a Pfaffian orientation.*

*Proof.* We claim that each $G_i$ is contained in $G$. It suffices to prove this for $i = 1$. Let the vertices of $C$ be $u_1, u_2, u_3, u_4$ in order. By 8.2 $G_2$ has two vertex-disjoint paths $P_1, P_3$ with ends $u_1, u_2$ and $u_3, u_4$, respectively, such that $G_2 \backslash V(P_1 \cup P_3)$ has a perfect matching. Similarly, $G_3$ has such paths $P_2, P_4$ with ends $u_2, u_3$ and $u_1, u_4$, respectively. Then $(G_1 \backslash E(C)) \cup P_1 \cup P_2 \cup P_3 \cup P_4$ is a model of $G_1$ in $G$, as required. Thus each $G_i$ is contained in $G$, and hence has a Pfaffian orientation by 6.4. □

If $G$ is a graph, we say that a set $X \subseteq V(G)$ is a *trisector* in $G$ if $|X| = 4$ and $G \backslash X$ has at least three components. The next lemma explains the significance of trisectors.

8.4. *Let $G$ be a connected brace, and let $X \subseteq V(G)$ with $|X| = 4$. Then $X$ is a trisector in $G$ if and only if $G$ can be expressed as a trisum of $G_1$, $G_2$ and $G_3$ along $C$, where $V(C) = X$.*

*Proof.* The "if" part follows immediately. For the "only if" part let $X$ be a trisector. It suffices to notice that 3.1 implies that $G \backslash X$ has a perfect matching, and hence $X$ contains two vertices from each color class. □

Theorems 6.5, 6.8, 8.2 and 8.4 imply an $O(|V(G)|^5)$ algorithm to test if a bipartite graph $G$ has a Pfaffian orientation. To improve its running time we need a few lemmas about trisectors.

8.5. *Let $G$ be a connected brace, and let $X$ be a set of vertices of $G$. Assume that $|N_G(X)| \leq 3$ and that $V(G) - X - N_G(X)$ contains vertices of both color classes of $G$. Then $|X| \leq 1$.*

*Proof.* Let $(A, B)$ be a bipartition of $G$. We may assume that $X \cap A \neq \emptyset$; then by 3.1 and the inclusion $N_G(X \cap A) \subseteq (X \cup N_G(X)) \cap B$,

(∗) $\qquad |X \cap A| + 2 \leq |N_G(X \cap A)| \leq |X \cap B| + |N_G(X) \cap B|.$

Similarly, if $X \cap B \neq \emptyset$, then $|X \cap B| + 2 \leq |X \cap A| + |N_G(X) \cap A|$, and it follows that $|N_G(X)| \geq 4$, a contradiction. Thus $X \cap B = \emptyset$, and (∗) implies that $|X| = |X \cap A| \leq |N_G(X)| - 2 \leq 1$, as required. □

8.6. *Let $G$ be a connected brace that has a Pfaffian orientation, and let $X, Y$ be trisectors in $G$. Then there exists a component $J$ of $G \backslash X$ such that $Y \subseteq V(J) \cup X$.*

*Proof.* Suppose for a contradiction that there exist two distinct components $J_1$ and $J_2$ of $G \backslash X$ and vertices $y_1 \in V(J_1) \cap Y$ and $y_2 \in V(J_2) \cap Y$.

We first show that $|V(J_i) \cap Y| = 1$ for $i = 1$ or $i = 2$. To this end suppose for a contradiction that $|V(J_i) \cap Y| > 1$ for $i = 1, 2$; then both those sets have cardinality two. By 8.2 there exist, for every component $L$ of $G\backslash Y$, two vertex-disjoint paths between $V(J_1) \cap Y$ and $V(J_2) \cap Y$ such that their vertex-sets except for their ends are contained in $V(L)$. Each of these paths intersects $X$, which is impossible, because $|X| = 4$ and yet $G\backslash Y$ has at least three components. This proves our claim that $|V(J_i) \cap Y| = 1$ for $i = 1$ or 2.

We may therefore assume that $V(J_1) \cap Y = \{y_1\}$. Let $L_1, L_2, \ldots, L_l$ be the components of $G\backslash Y$; then $l \geq 3$. By 8.2, for each $i = 1, 2, \ldots, l$, $V(L_i) \cup \{y_1, y_2\}$ includes the vertex-set of a path with ends $y_1$ and $y_2$. Thus each $V(L_i)$ intersects $X$, and hence $l \leq 4$, $|X \cap Y| \leq 1$, and $|X \cap V(L_i)| \leq 2$ with equality possible only if $X \cap Y = \emptyset$. Since for $i = 1, 2, \ldots, l$

$$N_G(V(L_i \cap J_1)) \subseteq (Y \cap V(J_1)) \cup (X \cap Y) \cup (X \cap V(L_i)),$$

it follows that $|N_G(V(L_i \cap J_1))| \leq 3$, and hence $|V(L_i \cap J_1)| \leq 1$ by 8.5. Since $V(J_1) = \{y_1\} \cup \bigcup_{i=1}^{l} V(L_i \cap J_1)$ we deduce that $|V(J_1)| \leq 5$, and so $J_1$ is a circuit of length four by 8.2. On the other hand $J_1\backslash y_1$ has no edges, a contradiction. □

8.7. *Let $G$ be a brace that has a Pfaffian orientation. Then $G$ has no subgraph isomorphic to $K_{2,3}$.*

*Proof.* Suppose for a contradiction that $K$ is a subgraph of $G$ that is isomorphic to $K_{2,3}$, and that $D$ is a Pfaffian orientation of $G$. Each of the three circuits of $K$ are central in $G$, and yet it follows (by considering the three paths of $K$ joining the vertices of $K$ of degree three) that at least one of them is not oddly oriented in $D$, a contradiction. □

8.8. *Let $G$ be a connected brace that has a Pfaffian orientation, let $G$ be a trisum of $G_1$, $G_2$ and $G_3$ along $C$, and let $X \subseteq V(G)$ with $X \neq V(C)$. Then $X$ is a trisector in $G$ if and only if it is a trisector in $G_1$, $G_2$ or $G_3$.*

*Proof.* Let $X$ be a trisector in $G_1$, say. We claim that the vertices of $V(C) - X$ belong to the same component of $G_1\backslash X$. To prove this claim suppose to the contrary that some two vertices of $V(C) - X$ belong to different components of $G_1\backslash X$. Then $V(C) \cap X$ consists of two diagonally opposite vertices of $C$. Let $G'_1$ be obtained from $G_1$ by joining nonadjacent pairs of vertices in $X$ that belong to opposite color classes. Thus $X$ induces a circuit in $G'_1$, and it follows from 8.2 and 6.5(ii) that $G'_1$ has a Pfaffian orientation. On the other hand, let $u \in V(C) - X$. Then $X \cup \{u\}$ is the vertex-set of a subgraph of $G'_1$ isomorphic to $K_{2,3}$, contrary to 8.8. This proves our claim that the vertices of $V(C) - X$ belong to the same component of $G_1\backslash X$. It follows that $X$ is a trisector of $G$, as desired.

Conversely, suppose that $X$ is a trisector of $G$. Then $X \subseteq V(G_i)$ for some $i \in \{1,2,3\}$ by 8.6. By 8.2 the vertices of $V(C) - X$ belong to the same component of $G \backslash X$, and hence $X$ is a trisector in $G_i$, as required. □

8.9. *Let $G$ be a connected brace on $n \geq 5$ vertices that has a Pfaffian orientation. Then $G$ has at most $n - 5$ trisectors.*

*Proof.* . We proceed by induction on $n$. If $G$ has no trisectors, then the result holds, and so we may assume that $G$ is a trisum of $G_1$, $G_2$ and $G_3$. By 8.2 each of $G_1$, $G_2$ and $G_3$ has a Pfaffian orientation, and hence $G_i$ has at most $|V(G_i)| - 5$ trisectors. By 8.9 $G$ has at most

$$|V(G_1)| - 5 + |V(G_2)| - 5 + |V(G_3)| - 5 + 1 \leq n - 5$$

trisectors, as required. □

## 9. An algorithm for Pfaffian orientations

We begin with the following easy algorithm.

9.1. ALGORITHM.
*Input. A bipartite graph $G$ with $m$ edges and a perfect matching $M$ of $G$.*
*Output. The set of all edges of $G$ that belong to no perfect matching of $G$.*
*Running time. $O(m)$.*

*Description.* Let $(A, B)$ be a bipartition of $G$, and let $D$ be the directed graph obtained from $G$ by directing every edge from $A$ to $B$, and contracting every edge in $M$. The problem is equivalent to finding the strongly connected components of $D$, which is well-known, and is described, for instance, in [2].

9.2. ALGORITHM.
*Input. A connected brace $G$ on $n$ vertices, and a list $\mathcal{L}$ of all trisectors of $G$.*
*Output. Either a Pfaffian orientation of $G$, or a valid statement that $G$ has no Pfaffian orientation.*
*Running time. $O(n^3)$.*

*Description.* If $|\mathcal{L}| > n - 5$, then we output the statement that $G$ has no Pfaffian orientation, and stop. By 8.9 this statement is correct. We now assume that $\mathcal{L}$ is empty. In this case we use a linear planarity algorithm such as [26] to either find a planar drawing of $G$, or determine that $G$ is nonplanar. If we find a planar drawing of $G$ we use Kasteleyn's algorithm [6], [7] (see also [13]) to output a Pfaffian orientation of $G$ and stop. If $G$ is nonplanar, then we check if $G$ is isomorphic to the Heawood graph. If it is, then we output a Pfaffian orientation of $G$ as in 6.3. If $G$ is not isomorphic to the Heawood graph, then by 6.8 and 8.4 it has no Pfaffian orientation. We output that information and stop. This completes the case when $\mathcal{L}$ is empty.

Thus we may pick a trisector $X \in \mathcal{L}$, and, by 8.4, express $G$ as a trisum of $G_1$, $G_2$ and $G_3$ along $C$, where $V(C) = X$. By 8.8 every member of $\mathcal{L} - \{X\}$ is a trisector in one of $G_1, G_2, G_3$, and all trisectors of $G_1, G_2, G_3$ belong to $\mathcal{L}$. Further, $G_1, G_2, G_3$ are braces by 6.5(iii), and $G$ has a Pfaffian orientation if and only if each of $G_1, G_2, G_3$ does by 6.5(ii) and 8.2. We apply 9.2 to each $G_1, G_2, G_3$ and an appropriate subset of $\mathcal{L}$. If each $G_i$ has a Pfaffian orientation, then we combine them as in the proof of 6.5(ii) to yield a Pfaffian orientation of $G$ and stop. If one of $G_1, G_2, G_3$ does not have a Pfaffian orientation, then neither does $G$. We output that information and stop.

9.3. ALGORITHM.
*Input. A connected brace $G$ on $n$ vertices.*
*Output. Either a Pfaffian orientation of $G$, or a valid statement that $G$ has no Pfaffian orientation.*
*Running time. $O(n^3)$.*

*Description.* If $G$ has more than $2n - 4$ edges, then it does not have a Pfaffian orientation by 7.3. We output that information and stop. Thus we may assume that $G$ has at most $2n - 4$ edges. For every $u, v \in V(G)$ we use the algorithm of Hopcroft and Tarjan [5] to find all trisectors $X$ of $G$ with $u, v \in X$. Thus we find all trisectors of $G$ in time $O(n^3)$, and apply 9.2.

Let $G$ be a connected 1-extendable bipartite graph. Let $\mathcal{C}_0 = \{G\}$, and assume that the sets $\mathcal{C}_0, \mathcal{C}_1, \ldots, \mathcal{C}_{i-1}$ of 1-extendable graphs have already been defined. If every member of $\mathcal{C}_{i-1}$ is a brace we stop; otherwise we choose $H \in \mathcal{C}_{i-1}$ that is not a brace. Then $H$ can be expressed as a 2-sum of two smaller connected 1-extendable bipartite graphs $H_1$ and $H_2$, and we put $\mathcal{C}_i = (\mathcal{C}_{i-1} - \{H\}) \cup \{H_1, H_2\}$. Let $k$ be the integer such that this process terminates with $\mathcal{C}_0, \mathcal{C}_1, \ldots, \mathcal{C}_k$. We say that $\mathcal{C}_k$ is a *decomposition of $G$ into braces*. Lovász [12] has shown that $\mathcal{C}_k$ is independent of the order in which individual graphs are decomposed, but we will not need that here. All we need is the following, which follows from 7.1.

9.4. *Let $G$ be a 1-extendable bipartite graph, and let $\mathcal{C}$ be a decomposition of $G$ into braces. Then $G$ has a Pfaffian orientation if and only if each member of $\mathcal{C}$ does.*

Let $G$ be a 1-extendable bipartite graph, and let $e$ be an edge of $G$ with ends $u_1, u_2$. We say that $e$ is *reducing* if $G\backslash u_1\backslash u_2$ is not 1-extendable. If $G_1$, $G_2$ and $e$ are as in the definition of 2-sum, we say that $G$ is a *2-sum of $G_1$ and $G_2$ along $e$*.

9.5. *Let $G$ be a 1-extendable bipartite graph, and let $M$ be a perfect matching of $G$.*

(i) *If $G$ is not a brace, then $M$ contains a reducing edge.*
(ii) *If $G$ is a 2-sum of $G_1$ and $G_2$ along $f \in M$ and $e \in M - \{f\}$ is reducing in $G_1$, then $e$ is reducing in $G$.*

*Proof.* To prove (i) let $e, Y_1, Y_2, u_1$ and $u_2$ be as in the definition of 2-sum, and let $F$ be the set of all edges with one end in $Y_1 \cup \{u_2\}$ and the other in $Y_2 \cup \{u_1\}$. Then $M \cap F \neq \emptyset$, and every edge of $F$ is reducing. Thus (i) follows.

To prove (ii) let $e$ have ends $u_1, u_2$, and let $(A, B)$ be a bipartition of $G$. Since $e$ is reducing in $G_1$, there exists a nonempty proper subset $X$ of $A \cap V(G_1) - \{u_1\}$ such that $|Y_1| = |X|$, where $Y_1 = N_{G_1}(X) - \{u_2\}$. By replacing $X$ by $B \cap V(G_1) - Y_1 - \{u_2\}$ we may assume that no end of $f$ belongs to $X$. Then $Y_1 = N_G(X) - \{u_2\}$, and hence $e$ is reducing in $G$, as desired. □

9.6. ALGORITHM.
*Input.* A connected 1-extendable bipartite graph $G$ with $2n$ vertices and $m$ edges, and a perfect matching $M$ of $G$.
*Output.* A decomposition of $G$ into braces.
*Running time.* $O(nm)$.

*Description.* Let $\mathcal{C}_0 = \{G\}$, and let $M = \{e_1, e_2, \ldots, e_n\}$. We repeat the following step for $i = 1, 2, \ldots, n$ and for every $H \in \mathcal{C}_{i-1}$. If $e_i \notin E(H)$ we put $H$ into $\mathcal{C}_i$; otherwise letting $u_1$ and $u_2$ be the ends of $e_i$ we find all subsets $X$ as in the definition of 2-sum. This is equivalent to finding the strongly connected components of $D(H\backslash u_1\backslash u_2, E(H) \cap M - \{e_i\})$, and hence can be done in time $O(m)$, as pointed out earlier. Using $\oplus$ to denote 2-sum along $e_i$, this gives a way to express $H$ as $(\cdots((H_1 \oplus H_2) \oplus H_3) \oplus \cdots) \oplus H_k$ for some graphs $H_1, H_2, \ldots, H_k$ such that for all $j = 1, 2, \ldots, k$, the graph $H_j\backslash u_1\backslash u_2$ is 1-extendable. We put $H_1, H_2, \ldots, H_k$ into $\mathcal{C}_i$.

After the last iteration the set $\mathcal{C}_n$ is a desired decomposition by 9.5.

9.7. ALGORITHM.
*Input.* A bipartite graph $G$ on $n$ vertices.
*Output.* Either a Pfaffian orientation of $G$, or a valid statement that $G$ has no Pfaffian orientation.
*Running time.* $O(n^3)$.

*Description.* We use the algorithm of Hopcroft and Karp [4] to find a perfect matching $M$ in $G$. If $G$ has no perfect matching, then every orientation of $G$ is Pfaffian. In that case we output an arbitrary orientation of $G$ and stop. Otherwise we use 9.1 to delete edges of $G$ that belong to no perfect matching of $G$. (Those edges may be directed arbitrarily in a Pfaffian orientation of the original graph.)

Now $G$ is 1-extendable. For every component $H$ of $G$ we proceed as follows. We find a decomposition $\mathcal{C}$ of $H$ into braces using 9.6, and apply 9.3 to every member of $\mathcal{C}$; if every member of $\mathcal{C}$ has a Pfaffian orientation, then those can be combined to give a Pfaffian orientation of $H$. If some member of $\mathcal{C}$ has no Pfaffian orientation, then neither does $H$.

If some component $H$ of $G$ has no Pfaffian orientation, then neither does $G$; otherwise the Pfaffian orientations of components of $G$ can be combined to yield a Pfaffian orientation of $G$.

*Acknowledgements.* We proved the main theorem of this paper in the summer of 1996. In the fall of that year, William McCuaig informed us that he had already obtained the same result, but had not publicly announced it. Consequently, a joint announcement of the result was made in [16]; but since the proof methods are quite different, we agreed with McCuaig that both papers should be submitted for publication independently.

We would like to thank an anonymous referee for carefully reading the manuscript, and for several helpful suggestions.

Ohio State University, Columbus, OH
*E-mail address*: robertso@math.ohio-state.edu

Bellcore, Morristown, NJ
*Current address*: Princeton University, Princeton, NJ
*E-mail address*: pds@math.princeton.edu

Georgia Institute of Technology, Atlanta, GA
*E-mail address*: thomas@math.gatech.edu